\documentclass[12pt,b5paper]{book}
\usepackage{lexicon}

\makeindex

\title{
\textbf{Arithmetical Foundations}\\
Part I Recursion\\
\begin{Huge}$\LOGO$\ \ $\mbf{A}$\end{Huge}
}

\author{\\
\qquad\\
\quad 
\\ 
\\  
M.\ Pfender
}

\date{
November 14, 2026}

\begin{document}
\maketitle

\noindent  
PD Dr.\ Michael Pfender\\
DMV, member no.\ 225\\
michael.pfender@gmx.de\\
alumnus TU Berlin

\bigskip\noindent   
With the cooperation of

\smallskip\noindent
stud.\ math.\ Anton Ernst\\
TU Berlin

\smallskip\noindent
Robert Sauer-Ernst\\
manager of the author

\smallskip\noindent
Dr.\ Wolfgang Wenning\\
alumnus TU Berlin

\vfill

\noindent 
The following students have contributed seminar talks:
           Sandra Andrasek,
           Florian Blatt,
           Nick Bremer,
           Alistair Cloete,
           Joseph Helfer,
           Frank Herrmann,  
           Julia Jonczyk,
           Sophia Lee,
           Dariusz Lesniowski,
           Wiktor Matysiak,
           Gregor Myrach,
           Thomas Richter,
           Olivia R\"ohrig,
           Paul Vater, and
           Jörg Wleczyk.

\vfill

\noindent

\paragraph*{abstract.}We fix constructive foundations for Arithmetic
on a map theoretical, algorithmic level.
In contrast to elementhood and quantification based 
traditional foundations such as Principia Mathematica $\PM,$ 
Zermelo-Fraenkel \textbf{set} theory $\ZF,$ or v.\,Neumann-Gödel-Bernays
\textbf{set} theory $\mbf{NGB},$
our fundamental primitive recursive theory $\PR$ has as its 
basic ``undefined" (not further defined) terms just terms for objects 
and maps. Goodstein's axiomatic uniqueness rules are
derived as theorems in $\PR$ -- with equality definability as a
side result --, $\bfS$ is obtained by predicate comprehension,
and $\hatS\iso\muS$ identifies $\mu$-recursive maps and while-loops
with partially defined maps.

\bigskip


\tableofcontents

\chapter*{Introduction}
\addcontentsline{toc}{chapter}{Introduction}
\markboth{\ \hrulefill\ Introduction}{Introduction\ \hrulefill\ }

We fix constructive foundations for Arithmetic
on a map theoretical, algorithmic level.
In contrast to elementhood and quantification based 
traditional foundations such as Principia Mathematica $\PM$\,, 
Zermelo-Fraenkel \textbf{set} theory $\ZF,$ or v.\,Neumann-Gödel-Bernays
\textbf{set} theory $\mbf{NGB}$\,,
our fundamental primitive recursive theory $\PR$ has as its 
basic ``undefined" (not further defined) terms just terms for objects 
and maps. 
On that language level it is \emph{variable free,} and it is free from
formal quantification over individuals like numbers or number pairs.
 
For an initial understanding of this formally variable-free 
\emph{free-variables} cartesian category theory approach, you may
consult the appendix to chapter 8 on \emph{Free variables as projections}.

Theory $\PR$ is strongly finitistic without formal
quantification, in the sense of 
\NAME{Skolem} 1919/1970, p.\ 153.\footnote{
  Was ich nun in dieser Abhandlung zu zeigen wünsche ist folgendes:
  
  \emph{Faßt man die allgemeinen Sätze der Arithmetik als    
  Funktionalbehauptungen}
  
  \emph{auf, und basiert man sich auf der rekurrierenden Denkweise, 
  so läßt sich}
  
  \emph{diese Wissenschaft in folgerichtiger Weise ohne Anwendung 
  der Russell-}
  
  \emph{Whitehead'schen Begriffe ``always'' und ``sometimes'' begründen}
  }

$\PR$ is a formal, combinatorial\footnote{
``combinatorial" in the sense of a game with numbers, signs, letters, symbols}
 category 
with cartesian \ie\ universal\index{universal} {product} and a natural numbers 
object (NNO\index{NNO}) $\N,$ a \emph{p.r.\ cartesian category,} 
cf.\ \NAME{Rom\`an}\index{Roman} 1989.

The NNO\index{NNO} $\N$ admits iteration of endo maps and the 
\emph{full schema of primitive recursion.} Such NNO\index{NNO} has been 
introduced in categorical terms by \NAME{Freyd}\index{Freyd} 1972, 
on the basis 
of the NNO\index{NNO} of \NAME{Lawvere}\index{Lawvere} 1964, and named later 
(\eg\ by \NAME{Maietti\index{Maietti}} 2010) \emph{parametrised NNO\index{NNO}.}

We remain on the purely {syntactical} level of this 
categorical theory and later {extensions:} 
{no formal semantics} necessary into an outside, non-combinatorial 
world, cf.\ Hilbert's formalistic program.\footnote{
  \NAME{Hilbert}\index{Hilbert} 1900/1970.
  }

Fundamental (categorical) \pr\ theory $\PR$ is developed from the endomap 
iteration schema $(\S)$ of \NAME{Eilenberg, Elgot}\index{Eilenberg, Elgot} 
1970.
We take as additional {axiom} \NAME{Freyd}'s\index{Freyd}
uniqueness of the initialised iterated endo map. This
gives the full schema of primitive recursion up to 
{uniqueness} of \pr\ maps defined by that schema.

Into our {variable-free} setting are introduced free variables, 
formally interpreted as {names} for identity and projection maps.
As a consequence, we have in the present context 
`{free variable}' as a defined notion. We have object
and map {constants} such as {terminal object,}
NNO\index{NNO}, zero constant and successor map, and use free metavariables for 
objects and for maps. 

Fundamental arithmetic is further developed along 
\NAME{Goodstein}'s\index{Goodstein} 1971 \textsc{free variables arithmetic} whose 
axiomatic uniqueness rules are
derived as theorems in $\PR$\,;
this gives the expected structure theorem for algebra and order
on NNO\index{NNO} $\N.$ ``On the way'', via Goodstein\index{Goodstein}'s truncated subtraction 
and his symmetry of the maximum function, we obtain the
\emph{equality definability theorem:} If predicative equality of
two p.r.\ maps is derivably true (``on all arguments''), then (objective) map equality \emph{between} 
these maps is derivable. 

The game is enriched by an (embedding, hence conservative\footnote{
  a conservative theory extension admits no new theorems in terms of the                
  extended theory.
  }) 
extension of theory $\PR$ by comprehension of 
predicates into  new (sub) objects. 

This enrichment makes emerging theory 
$\bfS=\PR+(\compr)$ with predicate comprehension more comfortable, 
in direction 
to \textbf{set} theory with its sets and subsets. 
 
$\bfS$ has a two-element 
boolean algebra $\two\iso\set{0,1}\subset\N$\,. 

Over theory $\bfS$
is constructed a theory $\hatS$ of \emph{partial} p.r.\footnote{
  primitive recursive  
  }\ maps\footnote{
  special Korrespondenzen in \NAME{Brinkmann, Puppe}
  \index{Brinkmann, Puppe} 1969
  } 
with half-terminal diagonal symmetric monoidal structure
in the sense of \NAME{Budach, Hoehncke} 1975;
$\mu$-\emph{recursive maps} and $\ttt{while}$ loop \ttt{program}s 
turn out to be just partial p.r.\  maps.

Finally we discuss, within $\hatS$\,, \emph{totally defined} maps i.e.\
partial maps with retractive domain of definition enumeration:
\textbf{Totality Lemma}.  

\part{RECURSION}

\chapter{Cartesian language}
\markboth{\ \hrulefill\ Cartesian language}{Cartesian language\ \hrulefill\ }

We develop from scratch the free-variables categorical
language with
terminal object\index{terminal object} $\one$   
and the data
of Natural Numbers Object NNO\index{NNO} $\N$ 
which comes with 
zero map\index{zero map} $0:\one\to\N$ and successor
\index{successor} endo map\index{endo map} 
$\s:\N\to\N$\,.

A special r\^ole is played by object $\one$\,,
in particular source of pointer arrows
\index{pointer arrow}. It works
as the domain object\index{domain object} for concrete 
``elements\index{element}'', ``constants\index{constant}'' 
$\bs{a}:\one\to A$
of any object $A$\,, in particular for 
concrete numbers\index{concrete numbers} 
 -- number constants\index{number constant} -- 
$\bs{n}:\one\to\N$\,. 

The set of fundamental symbols\index{fundamental symbol} of ``our'' 
categorical\index{categorical} language\index{language} is  
\begin{align*}
& \set{\id,\circ,\one,\N,\0,\s,\Pi},
\ \text{and equality sign `$=$'}.
\end{align*}

Identity\index{identity} $\id$ is the family of 
identity maps\index{identity map} for all objects,
$\circ$ is map 
composition\index{composition}, occasionally replaced by 
concatenation\index{concatenation},

$\Pi$ symbolises the family of terminal maps
\index{terminal map}
to object $\one$\,,

and `$=$' is to be equality\index{equality} of (\emph{between}) 
maps.

We give here the axioms\index{axiom} of (our) category 
theory\index{theory} in a fully 

formal way using Gentzen bars\index{Gentzen bar} for 
expression of inferences\index{inference}. 

\inference{[\,\one,\ \N\,]}
{(empty antecedent\index{antecedent} for this inference\index{inference})}
{$\mr{Obj}\ \one,\ \N$\\ 
&one-element object\index{one-element object}
  and\\
&natural number\index{natural number}s object NNO\index{NNO}\\
&$\mr{map}\ \0:\one\to\N\ \mr{zero} 
\ \mr{constant}$\\ 
&$\mr{map}\ \s:\N\to\N$ successor \index{successor} function
}\footnote{
  $\one$ is the terminal object\index{terminal object}, in the partial-
  maps context: a half-terminal\index{half-terminal} object 
  $\one=\set{\0}\subset\N$, see \textbf{chapter 6} below
  }

\inference{[\,\id\,]}
{$\mr{Obj}\ A$}
{$\mr{map}\ \id_A=\id:A\to A$\\ 
&identity\index{identity} map
}

\inference{[\emph{reflexivity}\index{reflexivity}]}
{$\mr{map}\ f$}
{$f=f$}

\inference{[\emph{symmetry}\index{symmetry}]}
{$\mr{map}\ f,g$\\
&$f=g$
}
{$g=f$}

\inference{[\emph{transitivity}\index{transitivity}]}
{$\mr{map}\ f, g, h$\\
&$f=g;\ g=h$
}
{$f=h$}

\inference{[\circ]}
{$f:A\to B;\ g:B\to C$}
{$\mr{map}\ (g\,f)=(g\circ f)=g(f):A\to C$\\
&map composition\index{composition}\\
&(outmost brackets may be omitted)
}

\inference{[\circ\,\mr{sub]}}
{$f, f':A\to B;\ g:B\to C;$ $f=f'$}
{$g\circ f = g\circ f'$ Leibniz'\index{Leibniz} 
substitutivity\index{substitutivity}:\\
&Substitution\index{substitution} of equals into same gives equals
}

\inference{[\mr{sub}\,\circ]}
{$f:A\to B;\ g, g':B\to C;$ $g=g'$}
{$g\circ f=g'\circ f$ second Leibniz'\index{Leibniz} 
substitutivity\index{substitutivity}:\\
&Substitution\index{substitution} of same into equals gives equals
}

\inference{[\circ\ \id]}
{$f:A\to B$}
{$f\circ\id=f\circ\id_A=f;$\\
&$\id\circ f=\id_B\circ f=f$\\
&neutrality\index{neutrality} of identities\index{identity} to composition
\index{composition}
}

\medskip
It follows a first statement on the use of (free) 
variables\index{free variable}.

\inference{\mbf{var}\ \mr{use}}
{$f:A\to B$\\
&$\mbf{var}\ a\,\mbf{\in}\, A,\ a:=\id_A$
}
{$f(a) = f(\id_A) \bydefeq f\circ\id_A=f$\\
&free variable\index{free variable} as identity\index{identity},\\
&$f(a)\,\mbf{\in}\, B$ ``dependent variable\index{dependent variable}''
}

\medskip 
Next axiom\index{axiom} is \emph{associativity\index{associativity} of composition\index{composition}.}

\medskip 
\inference{[\mr{ass}\ \circ]}
{$f:A\to B;\ g:B\to C;\ h:C\to D$}
{$(h\circ g)\circ f=h\circ(g\circ f):A\to D$\\
&$ = h\circ g\circ f = h\,g\,f = h(g(f(a)))$
}

\paragraph{Counting Remark.} Up to insertion of 
(composition-neutral) 
identities\index{identity}, the maps\index{map} of 
category theory\index{category theory} generated over 
$\s:\N\to\N$ are just the iterated maps\index{iterated map} 
$\s\circ\ldots\circ\mr{s}:\N\xto{\s}\ldots\xto{\s}\N$
of the successor\index{successor} map, as well as the
``external'' numerals\index{numeral} -- special number
constants\index{number constant}:

\inference{}
{(empty antecedent\index{antecedent})}
{$\0:\one\to\N$ numeral\index{numeral}}

\bigskip 
\inference{}
{$\bs{n}:\one\to\N$ numeral\index{numeral}}
{$(\s\circ\bs{n}):\one\to\N$ numeral\\
&\textbf{example}: $\ttt{3}=(\s\circ(\s\circ(\s\circ\0)))$
}

\chapter{Primitive Recursion}
\markboth{\ \hrulefill\ Primitive Recursion}{Primitive Recursion\ \hrulefill\ }

We introduce Gödel's\index{Gödel} primitive recursion -- called by him just
recursion\index{recursion}\footnote{later Ackermann found a \emph{recursive}
  function which is \emph{not} primitive recursive
}
--, beginning with the iteration schema in Eilenberg, Elgot\index{Eilenberg, Elgot} 1970.
We show the full schema of primitive recursion and
uniqueness of the NNO\index{NNO} $\N$ within the categorical theory $\PR$
of primitive recursion to be described in this chapter.

Skolem's\index{Skolem} influential introduction of recursive
arithmetic motivated R.\ L.\ Goodstein\index{Goodstein} to define
an entirely logic-free version of primitive recursive
arithmetic which we take up here by free-variable
Primitive Recursive Cartesian Category Theory $\PR$
with NNO\index{NNO}. We then extend $\PR$ by 
predicate-into-subobject comprehension to theory $\bfS$
where a boolean truth object $\two$ is constructed
together with a categorical free-variable approach
which allows to derive theorems of number theory
categorically.
We extend $\bfS$ into the categorical
theory $\hatS$ of \emph{partial 
recursive maps}, and introduce $\mu$-recursion and 
while-loops categorically.

Two-valued logic is based on a truth object\index{truth object} $\two$
with truth values $\false, \true:\one\to\two$ which together  
may be represented by the \textsc{diagram} 

\medskip
\quad $\xymatrix{\one\ar[r]^\false&\two&\one\ar[l]_\true}$\footnote{
  introduced formally as a \emph{weak coproduct} diagram in the
  theory $\bfS=\PR+(\mr{compr})$ of primitive recursion with
  comprehension of predicates into subobjects of $\PR$ objects,
  see \textbf{section 4.3} below
  }

\medskip
with boolean operations\index{boolean operations} shown in the 
\textsc{diagram}
\begin{align*}
&\xymatrix @+0em{
\one
\ar@<+0.5ex>[r]^\bot
\ar @<-0.5ex> [r]_\top 
&\two 
\ar[dr]^(.5)\neg
&&\two\times\two
\ar@<-1.5ex>@/_1.5pc/[ll]_(.6){\bs{\vert}}
\ar@<-1.0ex>@/_1.0pc/[ll]^(.6)\sminus
\ar@<-0.5ex>[ll]_{\dot{\lor}}
\ar@<+0.5ex>[ll]^\land
\ar@<-0.5ex>[dl]_(.5)\impl
\ar@<+0.5ex>[dl]^\biimpl
\\
&&\two
}
\end{align*}
In order to construct this truth object\index{truth object} 
$\two$ 
categorically\index{categorically} we will make use of the NNO\index{NNO}
(Natural Numbers Object) known from F.\ William
Lawvere,
\begin{align*} 
&\xymatrix@+0em{\one\ar[r]^\0&\N\ar[r]^\s&\N}
\end{align*}
on which we consider specific arithmetical 
operations\index{arithmetical operations} shown in the \textsc{diagram}
\begin{align*}
\xymatrix@+1em{
\one
\ar @<+0.5ex> [r]^\0
\ar @<-0.5ex> [r]_1 
& \N
\ar[d]^{\,> \0}
&\N\times\N
\ar@<-0.6ex>@/_1pc/[l]_{\dmin}
\ar@<-0.5ex>[l]_{+}
\ar@<+0.5ex>[l]^{\mul}
\ar@<-0.5ex>[ld]_{\leq}
\ar@<+0.5ex>[ld]^{\doteq}
\\
&\N
}
\end{align*}
constituting  a linearly ordered unitary 
commutative semiring\index{semiring} 
\begin{align*}
&\N=\bfan{\N,\one,\ttt{1},+,\mul}
\end{align*}
with 
zero $\0$\,, truncated subtraction\index{truncated subtraction}
`$\dmin$', order `$\leq$', and
equality\index{equality} predicate `$\doteq$'.

This diagram may be translated to a 
\emph{protoboolean\index{protoboolean}} one:
\begin{align*}
\xymatrix @+1em{
\one
\ar @<+0.5ex> [r]^\bot
\ar @<-0.5ex> [r]_\top 
& \N
  \ar[d]^{\sg}
&\N\times\N
\ar@<-0.6ex>@/_1pc/[l]_{\sminus}
\ar@<-0.5ex>[l]_{\dot{\lor}}
\ar@<+0.5ex>[l]^{\land}
\ar@<-0.5ex>[ld]_{\impl}
\ar@<+0.5ex>[ld]^{\biimpl}
\\
&\N
}
\end{align*}
In 1923 Skolem\index{Skolem} showed that large parts of 
elementary Number
Theory are to be reached by what would 
later be called primitive recursion\index{recursion}.   
Indeed Skolem\index{Skolem} laid out
the groundwork\footnote{
  \NAME{Skolem} 1923: ``Now what I wish to show in the
  present work is the following: If we consider the
  general theorems of arithmetic to be functional
  assertions and take the recursive mode of thought
  as a basis, then that science can be founded in a
  rigorous way without Russell and Whitehead's 
  notions `always' and `sometimes'.''} for 
\emph{Recursive Number Theory\index{number theory}}.

Herein Skolem\index{Skolem} uses the notation of Russell and Whitehead\index{Russell and Whitehead},
but restricts the scope of the quantifiers\index{quantifier} to finite domains
such that they loose their proper function of binding
variables\index{variable} in an unbounded way, and in principle become
eliminable by finite conjunction\index{conjunction}s 
and disjunction\index{disjunction}s.

Interestingly, in his concluding remarks, he sort of
excuses himself for retaining the logical quantifier\index{quantifier} 
notation -- if only restricted -- and announces a
formalisation ``free from this formal 
cumbrousness''\footnote{
  ``[On the other hand] I myself am no longer satisfied with 
  the way in which the logical elaboration proceeds here, 
  since by following the pattern of Russell and Whitehead's 
  work I made it too cumbersome from the purely formal point 
  of view. Yet, even in providing a foundations for 
  mathematics it is the substance that
  is important, not the notation. Soon I shall publish
  another work on the logical foundations of mathematics that 
  is free from this formal cumbrousness''
  },
i.e. the restricted quantifier\index{quantifier} notation. Skolem\index{Skolem} never
published what he announced, however he foresaw the
possibility of such an endeavor which was realised later
independently by R.\ L.\ Goodstein\index{Goodstein} 1941/45/57 and
H.\ B.\ Curry\index{Curry} 1941.

Goodstein\index{Goodstein} advanced this field in
``Recursive Number Theory: A development of Recursive
Arithmetic in a logic-free equation calculus''
(1957) in an important way which prompted Skolem\index{Skolem}
to review this book in the Journal of Symbolic Logic
\textbf{23}(2), pp.\ 227-228 (1958).

Skolem\index{Skolem} emphasises:
``[...] only the use of free variables remains, and this kind
of Arithmetic is called \emph{Recursive Number Theory}.''

Goodstein\index{Goodstein} introduces chapter 7 
\textsc{free-variable arithmetics}
of his 1971 book
in a programmatic way as follows:

``A substantial portion of arithmetic, including everything 
that would usually be called elementary number theory\index{number theory}, may be
codified in a calculus in which quantifiers\index{quantifier} play no part.
In view of the impossibility, which we established in the
previous chapter, of characterising existence\index{existence} and universality\index{universality}
in a formal system by means of quantifiers\index{quantifier}, quantifier-free
systems are of some importance, and this importance is
enhanced by the fact that such systems have, as we shall show,
a transparent freedom from contradiction\index{contradiction} which is not shared
by systems employing bound variables\index{bound variable}.
The codification of arithmetic which we shall describe is not
only entirely independent of predicate logic\index{predicate logic}, but it is
independent also of the axiom\index{axiom}s and inference rules of 
sentence logic\index{sentence logic}, and makes no appeal even to 
the axiom\index{axiom} of
induction\index{induction}. In our previous formalisation of arithmetic, 
arithmetic was based on formal logic\index{formal logic},
but in the new formalisation arithmetic is independent of
formal logic and sentence logic may, in turn, 
be based\footnote{
  in fact Goodstein constructs formal logic out of
  his logic-free recursive arithmetic including
  modus ponens. On this basis he then proves the
  schema of induction\index{induction} which we will rely on here.
  We give a proof of Peano induction\index{induction} 
  in section 5.4 Peano in $\bfS$} 
on arithmetic.''
  
This is what motivated us to develop Goodstein\index{Goodstein}'s approach 
in categorical terms.

Natural numbers made their way into 
category theory not before F.\ William Lawvere 1964
introduced his Natural Numbers Object\index{Natural Numbers Object}\footnote{
  \NAME{Lawvere} 1964, \textsc{Axiom\index{axiom} 3.}
  } 
NNO\index{NNO} by the \emph{simple recursion}

\textsc{diagram}
\begin{align*}
&\xymatrix@-.7em{
&&\N 
\ar[rr]^{\s}
\ar@{-->}[dd]^a
\ar@{}[ddrr]|{=}
&&\N 
\ar@{-->}[dd]^a& 
\\
\one
\ar[urr]^\0
\ar[drr]_{a_0}
\ar@{}[rrr]|(.4){=}
&&&&&(\mr{FWL})  
\\
&&A
\ar[rr]^t
&&A
}
\end{align*}

Here the sign `$=$' stands for diagram commutativity,
and the dashed arrows stand for the uniqueness\index{uniqueness} of the
\emph{sequence} 

$a:\N\to A$\,;\ $a_0=a\,\0,\ a_1=t\,a_0\,,\ a_2=t\,a_1$ etc., 

$a_0$ is the \emph{anchor\index{anchor}} of the \emph{recursion\index{recursion}}, 
$a_{\s\,n}=t\,a_n$\ \emph{refers back} to $a_n$\,.

\medskip\noindent
Inspired by 
\NAME{Eilenberg, Elgot}\index{Eilenberg, Elgot} 1970\footnote{   
  \emph{Recursiveness}, p.\ 13 ``right exponential''.
  We develop recursive Arithmetic categorically in the
  framework of \NAME{Goodstein} 1971 
  },
we introduce a schema $(\_ ^\S)$ of 
\emph{basic iteration}\index{basic iteration}, which yields a 
first form here of a parametrised NNO\index{parametrised NNO} by the 

\textsc{diagram}\footnote{
  in the late 1960ties S.\ Eilenberg turns from 
  Categories -- common work with S.\ MacLane -- to
  diagram-free Computer Science of J.\ B.\ Wright 
  and C.\ C.\ Elgot in the IBM Watson Research Center.
  The diagram -- with unique(!) commutative fill-in
  is our -- categorical -- contribution. Uniqueness is
  dropped in \NAME{Lambek} 1999 where
  a \emph{weak} Natural Numbers Object is considered
  }

\begin{align*}
\begin{split}
(\mr{EE})
\end{split}
\begin{split}
&\xymatrix@-.7em{
&A\times\N 
\ar[r]^{A\times\s}
\ar@{-->}[dd]^{k^\S}
\ar@{}[ddr]|{=}
&A\times\N 
\ar@{-->}[dd]^{k^\S}& 
\\
A
\ar[ur]^{(\id_A,\0)}
\ar[dr]_{\id_A}
\ar@{}[r]|{=}
&&&(\ ^\S)  
\\
&A
\ar[r]^k
&A
}
\end{split}
\end{align*}

In the special case where $k:A\to A$ is chosen to be the 
successor map\index{successor map} $\s:\N\to \N$\,, 
this schema gives in particular
addition recursively:
\begin{align*}
&a+n=+(a,n):=\s^\S(a,n)=\s^n(a):\N\times\N\to\N.\\
&a+\0=\s^\S(a,\0)=\s^\0(a)=a\ \,\text{and}\\
&a+\s(n)=\s^\S(a,\s n)=\s\,\s^\S(a,n)=\s(a+n)
\end{align*}

Above -- characteristic for our approach to 
categorical recursive foun\-dations -- the 
\emph{free variables}\index{free variable} $a$ and $n$\,
over $A$ and $\N$ respectively, are to
stand for the left projection\index{projection} 
$\ell:A\times\N\to A$ and the right projection\index{projection} 
$\err:A\times\N\to\N$ out
of the cartesian\index{cartesian} (universal\index{universal}) product
$A\times\N$\,.
  
A more general iteration\index{iteration} schema was 
introduced in Freyd 1977,
by replacing $A\xto{\id}A\xto{k}A$ in the above with
$A\xto{f}B\xto{g}B$\,.

Freyd \emph{initialised iterated\index{iterated} endomap\index{endomap} schema}
in free-variables notation:

\inference{\mr{(FR)_{FV}}}
{$f=f(a):A\to B$\ \,\emph{}\\
&$g=g(b):B\to B$\ \,\emph{iteration step}
}
{
$h=h(a,n):=g^\S(f(a),n):A\times\N\to B$\\
& \emph{initialised iterated}
}

with the defining properties
\begin{align*}
&h(a,\0)=f(a):A\to B\ \,\text{and}\\
&h(a,\s\,n)=g(h(a,n)):A\times\N\to B
\end{align*}

Variable-free\index{variable-free} notation of this schema
$(\mr{FR})_{\mr{FV}}$
is obtained by
replacing the free variables\index{free variables} $a$ and $n$ with the 
projections\index{projection}
$\ell:A\times\N\to A$ and
$\err:A\times\N\to\N$\ \,as follows.

\inference{\mr{(FR)}}
{$f:A\to B$\ \,\emph{}\\
&$g:B\to B$\ \,\emph{iteration step}
}
{
$h:=g^\S\circ(f\times\N):A\times\N\to B\to B$\\
&$= g^\S\circ(f\circ\ell_{A,\N},\err_{A,\N})$\\
&$:A\times\N\to B\times\N\to B$\ \,
\emph{initialised iterated\index{iterated}}
}

with defining properties
\begin{align*}
&h(\id_A,\0)=f:A\to B\ \,\text{and}\\
&h\circ(A\times\s)\\
&=h\circ(\ell_{A,\N},\s\circ\err_{A,\N})\\
&=g\circ h:A\times\N\to B\to B,
\end{align*}
by \textsc{diagram}:
\begin{align*}
&\xymatrix@-.7em{
& A\times\N
\ar[rr]^{A\times\s}
\ar[dd]_h^{g^\S\circ(f\times\N)}
& & A\times\N
\ar[dd]_h^{g^\S\circ(f\times\N)}
\\
A
\ar[ur]^{(\id_A,\0)}
\ar[dr]_f
\ar@{}[r]|(0.5){=}
&\ar@{}[r]|(1.35){=}
&&&&(\mr{FR})
\\
& B
\ar[rr]^g
& & B
}
\end{align*}
Thus we can use two different notations for the same
function here. In the first one the free variables\index{free variable}
$a$ and $n$ occur explicitly, while in the second one,
the \emph{variable-free\index{variable-free}} notation, they are
\emph{implicitly} present behind the curtain of
projections\index{projection}. It is the second notation which 
corresponds to Skolem\index{Skolem}'s and Goodstein\index{Goodstein}'s emphasis
on \textsc{variable-free}\index{variable-free} formalisations.

In what we will call Skolem\index{Skolem} category 
$\bfS$\,,
which allows the definition of subobjects\index{subobject} as
predicates\index{predicate}, the boolean truth object\index{truth object} $\two$ will
become definable such that in $\bfS$ what was above called
a protoboolean\index{protoboolean} diagram finally gets turned into a genuine 
boolean\index{boolean} one.

As is specified in the table of Contents, the first chapters
cover categorical\index{categorical} approaches to 
Arithmetic\index{Arithmetic}, and do not
transcend primitive recursion\index{recursion}. In the last three chapters
general recursion\index{general recursion} is introduced via partial functions\index{partial functions}
categorically\index{categorically}.

This section establishes the ``full''
schema $(\mr{pr})$ of Primitive Recursion within 
free-variables cartesian\index{cartesian} categorical\index{categorical} theory\index{theory}
$\PR$ of primitive recursion\index{recursion}. Schema $(\mr{pr})$ will
be constructed with the help of Freyd\index{Freyd}'s NNO\index{NNO}-schema
$(\mr{FR})$ and uniqueness\index{uniqueness} schema of the 
iterated\index{iterated} which
we call $(\mr{FR}!)$\,. 
All of this needs the free cartesian\index{cartesian}
category\index{category} generated over the NNO\index{NNO}-diagram
$\one\xto{\0}\N\xto{\s}\N$\,. Based on this we
consider the algebraic operations\index{algebraic operation} 
and the order\index{order}
structure on $\N$ in this section. We start with 
elementary definitions partly for reasons of 
self-containment, and also in order to fix notation\footnote{
  cf.\ \NAME{Eilenberg, Elgot} 1970 I,3 and I,4
  }.

\section{Cartesian Structure of $\PR$}

\paragraph*{Axiom.} For each object $A$ there is a
\emph{terminal} map\index{terminal map} (projection)
$\Pi_A:A\to\one$\,. 
$\Pi_A$ is to be unique\index{unique} as a map to $\one$\,.

\paragraph*{Axiom.} To any objects $A, B$ is associated
(recursively) a cartesian product object 
$A\times B$ with left and right projection\index{projection}
\begin{align*}
&\xymatrix{
A
&A\times B
\ar[l]_(.55)\ell
\ar[r]^(.55)\err
&B
}
\end{align*} 
 
which can be used to create variables\index{variable}

\inference{}
{$\mbf{Obj}\ A,B$\\
&$\mbf{var}\ a\ \mr{on}\ A,$\ $\mbf{var}\ b\ \mr{on}\ B$
}
{
$\mbf{map}\ a:=\ell_{A,B}:A\times B\to A$\\
&$\mbf{map}\ b:=\err_{A,B}:A\times B\to B$
}

\medskip 
\textsc{diagram}
\begin{align*}
&\xymatrix{
A
&A\times B
\ar[l]^(.55)\ell_(.55)a
\ar[r]_(.55)\err^(.55)b
&B
}
\end{align*}

\paragraph{Axiom.} \emph{Induced map\index{induced map} into a product}:

\inference{}
{$f:C\to A,\ g:C\to B$}
{$\mbf{map}\ (f,g):C\to A\times B$\\ 
&$\ell_{A,B}\circ(f,g)=f:C\to A$\\ 
&$\err_{A,B}\circ(f,g)=g:C\to B$
}

\medskip  
\textsc{diagram}
\begin{align*}
&\xymatrix@+0em{
\ar@{}[dr]|(0.6)= 
& A
\\
C
\ar @/^.5pc/[ru]^{f}
\ar[r]^(0.45){(f,g)} 
\ar @/_.5pc/ [rd]_{g}
& A \times B
\ar [u]_{\ell}
\ar [d]^{\err}
\ar@{}[dl]|(0.4)=
\\
& B
}
\end{align*}

\paragraph*{Axiom} \emph{Retractive pairing}:

\inference{}
{$h:C\to A\times B$}
{$h=(\ell_{A,B}\circ h,\err_{A,B}\circ h):C\to A\times B$ }

\medskip
Thus $h$ becomes an induced map\index{induced map}.

$$
\xymatrix{
\\
\ar[u]^y
\mul&\mul&\mul&\mul
&\mul
&
\\
2
&
\mul&\mul&(a,b)
\ar@{|->}[dd]^{\ell}_a
\ar@{|->}[lll]_{\err}^b
&\mul&
\\
\mul&\mul&\mul&\mul&\mul&
\\
\ar[uuuu]
\ar[rrrrr]
\mul&\mul&\mul&3
&\mul
\ar[r]^x
&
}
$$

\begin{center}
${\textbf{Cartesian product}}\ \N\times\N$

with projections

$a=\ell:\N\times\N\to\N$ and $b=\err:\N\times\N\to\N$

here $a(\ttt{3},\ttt{2})=\ttt{3}$\,, $b(\ttt{3},\ttt{2})=\ttt{2}$
\end{center}

\paragraph*{Retraction Remark.}\footnote{
  asked for by Mr Aimane Sarih} 
A \emph{retraction} $f:A\to B$ is a map
coming with a \emph{coretraction} $f':B\to A$ such that
\begin{align*}
&\xymatrix@+2em{
&&\\
B
\ar@/^2pc/[rr]^{\id_B}
\ar[r]_{f'}
&A
\ar[r]_f
\ar@{}[ur]|(.2)=
&B\\
b
\ar@{|->}[r]&
f'(b)
\ar@{|->}[r]&
f\circ f'(b)=\id_B(b)=b
}\\
&\text{$f'(b)$ is \emph{retracted} to $b$ by \emph{retraction} $f$}
\end{align*}
This is to explain the name ``retractive pairing'' for this \textbf{axiom}
on pairing: The pair $\ell_{A,B}\circ h$\,, $\err_{A,B}\circ h$
is \emph{retracted} by the \emph{induced\index{induced}} map 
formation $(\_\,,\_)$ to the
map $h:C\to A\times B$\,.

\subsubsection*{Theorem} 
\emph{Uniqueness}\index{uniqueness} of the induced map\index{induced map}.

\inference{}
{$f:C\to A,\ g:C\to B,\ h:C\to A\times B$\\
&$\ell_{A,B}\circ h=f:C\to A,\ 
\err_{A,B}\circ h=g:C\to B$}
{$h=(f,g):C\to A\times B$}
\paragraph*{Proof.}
\begin{align*}
& h:C\to A\times B\\ 
& = (\ell_{A,B}\circ h,\err_{A,B}\circ h):
C\xto{(\ell h,\err h)}A\times B\\ 
&\quad \text{by retractive pairing and then by antecedent}\\
&=(f,g):C\to A\times B\ \,\textbf{\qed}
\end{align*}

\paragraph*{Definition.} Binary cartesian product\index{cartesian product} of \emph{maps}:

\inference{}
{$f:A\to A',\ g:B\to B'$}
{$(f\times g):=(f\circ\ell_{A,B},g\circ\err_{A,B})
:A\times B\to A'\times B'$
}

\paragraph*{Definition.} Variables\index{variable} creation in a cartesian map product\index{cartesian product}:

\inference{}
{$f:A\to A',\ g:B\to B'$\\ 
&$\var\ a$ on $A$\,, $\var\ b$ on $B$
}
{$(f\times g)(a,b)=(f(a),g(b))
:=(f\circ\ell_{A,B},g\circ\err_{A,B})$\\
&$:A\times B\to A'\times B'$
}

\medskip
in particular
\begin{align*}
&a\times b=\id_A\times\id_B=(a,b)=(\ell_{A,B},\err_{A,B})
=\id_{A\times B}
\end{align*}

\textsc{diagram}
\begin{align*}
&\xymatrix@+0em{
A
\ar[r]^f
\ar@{}[dr]|=
& A'
\\
A\times B
\ar@{-->}[r]^{f\times g}
\ar[u]_\ell
\ar[d]^{\mr{r}}
\ar@{}[dr]|=
& A'\times B'
\ar[u]_\ell
\ar[d]^{\mr{r}}
\\
B
\ar[r]^g
& B'
}
\end{align*}

\subsubsection*{Theorem} 
\emph{Bifunctor}\index{bifunctor} property of the cartesian product\index{cartesian product}

\inference{}
{$f:A\to A',\ f':A'\to A'';$\\
&$g:B\to B',\ g':B'\to B''$
}
{$\id_A\times\id_B=\id_{A\times B}:A\times B\to A\times B$\\
&$(f'\circ f)\times(g'\circ g)=(f'\times g')\circ(f\times g):$\\
&$(A\times B)\to(A''\times B'')$
}
  
\paragraph*{Proof} by uniqueness\index{uniqueness} of the induced maps\index{induced map},

\medskip
\textsc{commutative diagram}
\begin{align*}
&\xymatrix@+0em{
A
\ar[r]^f
& A'
  \ar[r]^{f'}
& A''
\\
(A\times B)
\ar[u]_{\ell}
\ar[d]^{\err}
\ar[r]^{(f\times g)}
\ar@/_1pc/@{-->}[rr]_{((f'\,f)\times(g'\,g))}
& (A'\times B')
\ar[u]_{\ell}
  \ar[r]^{(f'\times g')}
\ar[d]^(.7){\err}
& (A''\times B'')
\ar[u]_{\ell}
\ar[d]^{\err}
\\
B
\ar[r]^g
& B'
  \ar[r]^{g'}
& B''
}
\end{align*}

\section{Freyd's initialised iteration}

\paragraph*{Axiom.} Lambek\index{Lambek} \emph{distributivity equation}.
\index{distributivity equation}
The composition operation `$\circ$' distributes over the
operation $(\_\,,\_)$ of forming the induced map
\index{induced map}:
 
\inference{}
{$h:D\to C,\ f:C\to A,\ g:C\to B$}
{$(f,g)\circ h=(f\circ h,g\circ h):D\to A\times B$}

\subsubsection*{Theorem}
With the help of basic (parametrised) iteration $k^\S$ 
(Eilenberg/Elgot's exponentiation) we get Freyd's
initialised iterated endomap schema $(\mr{FR})$
above.
\paragraph{Proof.}
Consider the following \textsc{diagram} wherein the above
``$A$-level'' is connected to the below ``$B$-level''
by the \emph{} function $f$\,.
\begin{align*}
&\xymatrix@+1em{
A
\ar[r]^{(\id_A,\0)}
\ar[d]_f
\ar@{-->}[dr]_(0.4){(f,\0)}
& A\times\N
\ar[r]^{A\times\s}
\ar[d]_(0.5){f\times\N}
\ar@{-->}[dr]^{f\times s}
\ar@{}[dl]|(0.2)=
\ar@{}[dl]|(0.75)=
& A\times\N
\ar[d]_{f\times\N}
\ar@{}[dl]|(0.25)=
\ar@{}[dl]|(0.65)=
\\
B
\ar[r]^{(\id_B,\0)}
& B\times\N
\ar[r]^{B\times\s}
& B\times\N
}
\end{align*}
The commutativity of the left square is established by 
the Lambek distributivity
\begin{align*}
&(\id_B,\0\circ\Pi_B)\circ f
=(\id_B\circ f,\0\circ\Pi_B\circ f)\\
&=(f,\0\circ\Pi_A):A\to B\times\N,
\end{align*} 
that of the right square by decomposition of the 
cartesian map product 
$f\times\s:A\times\N\to B\times\N$\,.

Replacing $k:A\to A$ in 
the basic iteration schema $(\mr{EE})=(\ ^\S)$
by the endomap $g:B\to B$\,, and 
in order to fit its gluing along
\begin{align*}
&B\xto{(\id_B,\0)}B\times\N\xto{B\times\s}B\times\N,
\end{align*}
supplemented by the
rounded arrows $f:A\to B$ and $h=g^\S\circ(f\times\N)$\,,
we get: 
\begin{align*} 
&\xymatrix@+1em{
A
\ar[r]^{(\id_A,\0)}
\ar[d]_f
\ar[dr]_(0.4){(f,\0)}
\ar@/_5pc/[ddr]_f
& A\times\N
\ar[r]^{A\times\s}
\ar[d]_(0.5){f\times\N}
\ar[dr]^{f\times s}
\ar@{}[dl]|(0.2)=
\ar@{}[dl]|(0.75)=
\ar@/^2pc/[dd]^(0.7)h
& A\times\N
\ar[d]_{f\times\N}
\ar@{}[dl]|(0.25)=
\ar@{}[dl]|(0.65)=
\ar@/^2pc/[dd]_(0.7)h^{g^\S\,\circ\,(f\times\N)}
\\
B
\ar[r]_{(\id_B,\0)}
\ar[dr]_{\id_B}
& B\times\N
\ar[r]_{B\times\s}
\ar[d]_{g^\S}
\ar@{}[dr]|(0.5)=
\ar@{}[dl]|(0.35)=
& B\times\N
\ar[d]_{g^\S}
\\
\ar@{}[ur]|(0.25)=
& B
\ar[r]^g
  & B
}
\end{align*}
The periphery of the \textsc{diagram} coincides with
$(\mr{FR})$\ \,\textbf{\qed}

\subsubsection*{Axiom. Uniqueness $(\mr{FR}!)$} 
The uniqueness\index{uniqueness} of the iterated endomap\index{endomap} $g$ initialised by
$f$ is postulated.

\inference{(\mr{FR!})}
{$f:A\to B,\ g:B\to B$\,,\\ 
&$h:A\times\N\to B$ arbitrary such that\\
&$h(a,\0)=f(a):A\to B$\ \,\text{and}\\
&$h(a,\s n)=(g\circ h)(a,n):A\times\N\to B\to B$
}
{$h(a,n)=g^\S(f(a),n)$\,.\\
&Elimination of the free variables turns this into\\
&$h=g^\S\circ(f\times\N):A\times\N\to B\times\N\to B$\,.
}

\subsubsection*{Theorem. Uniqueness of the basic iterated map}
From $(\mr{FR}!)$\,, by returning to 
$A\xto{\id_A}A\xto{k}A$ from $A\xto{f}B\xto{g}B$\,, 
uniqueness of the iterated $k^\S$ in the basic iteration 
schema $(\mr{EE})=(\ ^\S)$ is entailed.

\subsubsection*{Theorem. Recovering Lawvere's NNO}

Choosing $\one\xto{a_0}A\xto{t}A$ instead of 
$A\xto{f}B\xto{g}B$ in both $(\mr{FR})$ and $(\mr{FR}!)$
we get the Lawvere NNO\index{Lawvere NNO}, which we write here in equivalent form using the natural isomorphism
$A\iso\one\times A$\,.

\textsc{diagram}
\begin{align*}
&\xymatrix@-.7em{
&\one\times\N
\ar[rr]^{\one\times\s}
\ar@{-->}[dd]_a^{t^\S\circ(\id_\one\times\N)}
&&\one\times\N
\ar@{-->}[dd]_a^{t^\S\circ(\id_\one\times\N)}
\\
\one
\ar[ur]^{(\id_\one,\0)}
\ar[dr]_{a_0}
\ar@{}[r]|(0.5){=}
&\ar@{}[r]|(1.45){=}
&&&&
\\
&A
\ar[rr]^t
&&A
}
\end{align*}
\textbf{\qed}

The uniqueness\index{uniqueness} up to
natural isomorphism\index{natural isomorphism} of the NNO\index{NNO} is already inherent in
Lawvere's original schema $(\mr{FWL})$ of \emph{simple}
recursion\index{simple recursion}.

For the converse derivation, of $(\mr{FR})$ from
$(\mr{FWL})$\,, \NAME{Freyd} 1971 uses the higher order 
structure of his topoi\index{topos} with NNO\index{NNO}.

\section{The full schema (pr) of Primitive Recursion}

Schema (pr) is stated 
here first in free-variables notation.\footnote{
   cf.\ \NAME{Pfender, Kröplin, Pape} 1994}

\subsubsection*{Theorem} 
\emph{Recursion Schema $(\mr{pr})$}.

\inference{(\mr{pr})}
{$g=g(a):A\to B$\ \emph{},\\
&$h=h((a,n),b):(A\times \N)\times B\to B$
\ \emph{recursion step}
}
{$f(a,n)=[h^\star g](a,n):A\times\N\to B\ \text{such that}$\\
&Anchor: $f(a,\0)=g(a),$\\
&Step: $f(a,\s(n))=h((a,n),f(a,n)).$
}

\medskip
Here the general primitive recursive map $f$ is
anchored by $g$ and iteratively extended by Step $h$\,. This is 
expressed by
$[h^\star g]$\,, where the star stands for iteration of $h$\,.
 
This $f$ depends on $a$\ in\ $A$\,, the recursion 
parameter $n$ in $\N$\,, and the respective previously 
generated value\index{value} $b=f(a,n)$ in $B$\,.

\subsubsection*{Theorem} 
The following variable-free\index{variable-free} schema 
$(\mr{pr})$
is a consequence of schemata $(\mr{FR})$ 
and $(\mr{FR}!)$\,.\footnote{
  cf.\ \NAME{Freyd} 1972, proposition 5.22 -- here 
  brackets are added
  }

\inference {(\mr{pr})}
{$g:A \to B$\\
&$h:(A\times\N)\times B\to B$
}
{$f=[h^\star g]:A\times\N\to B$  given s.t.\\ 
&Anchor:
$f\circ(\id_A,\0)= g:A\to B$\\
&Step:
$f\circ(\id_A\times\s):A\times\N\to A\times\N\to B$\\
&$=h\circ(\id_{A\times\N},f)
:(A\times\N)\to(A\times\N)\times B\to B$
}

\paragraph*{Proof.}
That $f$ is to satisfy Anchor and Step 

is expressed by the \textsc{diagram}s Anchor and Step.
 
Anchor\qquad\qquad\qquad\qquad  Step
\begin{align*}
\begin{split} 
&
\xymatrix@-.5em{
& A\times\N 
  \ar[dd]^f
\\
A
\ar [ur]^{(\id,\0)}
\ar [dr]_g
\ar @{} [r] |{=}
&
\\
& B
} 
\end{split}
\begin{split}
&
\xymatrix@+1.8em{
A\times\N 
\ar[r]^{A\times\s}
\ar[d]_{(\id_{A\times\N},f)}
\ar@{}[dr]|{=}
&A\times\N 
\ar [d]^f
\\
(A\times\N)\times B
\ar[r]^(.5)h
&B
} 
\end{split}
\end{align*}

For the construction of such an $f:A\times\N\to B$ 

in Anchor and Step consider

the \textsc{diagram}
\begin{align*}
&
\xymatrix@-.5em{
& A\times\N 
\ar[rr]^{A\times\s}
\ar@{-->}[dd]^k
\ar@{}[ddrr]|=
&&A\times\N 
\ar@{-->}[dd]^k
\\
A
\ar[ur]^{(\id,\0)}
\ar[dr]^{\hat{g}}_{((\id,\0),g)}
\ar @{}[r]|(0.5)=
&&
\\
&(A\times\N)\times B
\ar [rr]^{\hat{h}}_{((A\times\s)\circ\ell,h)}
&&(A\times\N)\times B
}
\end{align*}
Here 
\begin{align*}
&\hat{g}:=((\id,\0),g):A\to(A\times\N)\times B\ \text{and}\\ 
&\hat{h}:=((A\times\s)\circ\ell,h):(A\times\N)\times B\to(A\times\N)\times B
\end{align*}
The dashed fill-in $k$ is obtained by (FR),
and $k$ is unique by axiom\index{axiom} schema (FR!).

Consider now the following extension of the previous

\textsc{diagram}:
\begin{align*}
&\xymatrix@+0em{
& A\times\N 
\ar[rr]^{A\times\s}
\ar[dd]_k
\ar@{}[ddrr]|=
\ar@{-->}@/^1.7pc/[ddd]^(.5){\ell\circ k}_{k_l}
&&A\times\N
\ar[dd]_k 
\ar@{-->}@/^1.7pc/[ddd]^(.5){\ell\circ k}_{k_l}
\\
A
\ar[ur]^{(\id,\0)}
\ar[dr]^(.5){\hat{g}}_(.6){((\id_A,\0),g)}
\ar@{}[r]|(0.5)=
\ar@{}[rr]|(.57)=
\ar@{}[rr]|(1.55)=
\ar@/_2pc/[ddr]_{(\id,\0)}
& &
\\
& (A\times\N)\times B
\ar[rr]^{\hat{h}}_{((A\times\s)\circ\ell,h)}
\ar[d]_{\ell}
\ar @{}[dl]|(0.3)=
&&(A\times\N)\times B  
\ar[d]_\ell
\\
& A\times\N
\ar[rr]^(.5){A\times\s}
&\ar@{}[u]|(.5)=
&A\times\N
}
\end{align*}

In the left part we have
\begin{align*}
&k\circ(\id_A,\0)=\hat{g}:A\to(A\times\N)\times B
\end{align*} 

by the properties of dashed $k$ in the previous diagram.

In the right part for the same reason holds 
\begin{align*}
&k\circ(A\times\s)
:A\times\N\to A\times\N\to(A\times\N)\times B\\
&=\hat{h}\circ k
:A\times\N\to(A\times\N)\times B\to(A\times\N)\times B
\end{align*}

For the dashed arrow on the left side of the diagram uniqueness  
$k_l:=\ell\circ k$ of $k$ comes from
\begin{align*}
&\ell\circ k\circ(\id_A,\0):A\to A\times\N\to (A\times\N)\times B\to A\times\N\\
&=\ell\circ\hat{g}=(\id_A,\0):A\to A\times\N.
\end{align*}

Correspondingly on the right side holds
\begin{align*}
&\ell\circ k\circ(A\times\s):A\times\N\to A\times\N\to(A\times\N)\times B\to A\times\N\\
&=\ell\circ\hat{h}\circ k\\
&=(A\times\s)\circ\ell\circ k:A\times\N\to(A\times\N)\times B\to A\times\N\to A\times\N.
\end{align*}

On both sides of the diagram what holds for $\ell\circ k$ also holds 

for $\id_{A\times\N}$\,. Therefore by (FR!) we get
\begin{align*}
&\ell\circ k=\id_{A\times\N}\qquad(*)
\end{align*}

Since the left projection $\ell$ is a retraction of $k$ 
we call $(*)$ the \emph{retraction property} of $\ell$ with 
respect to $k$\,.
 
In the following diagram the composition $\err\circ k$ of the right projection $\err$ with $k$ yields the desired
\begin{align*}
& f=\err\circ k:A\times\N\to(A\times\N)\times B\to B.
\end{align*}

To see this consider the \textsc{diagram}
\begin{align*}
&\xymatrix@+0em{
& A\times\N 
\ar[rr]^{A\times\s}
\ar[dd]_k
\ar@{}[ddrr]|=
\ar@/^1.7pc/[ddd]^(.5)f
&&A\times\N
\ar[dd]_k 
\ar@/^1.7pc/[ddd]^(.5)f
\\
A
\ar[ur]^{(\id,\0)}
\ar[dr]^(.5){\hat{g}}_(.6){((\id_A,\0),g)}
\ar@{}[r]|(0.5)=
\ar@{}[rr]|(.57)=
\ar@{}[rr]|(1.55)=
\ar@/_2pc/[ddr]_g
& &
\\
& (A\times\N)\times B
\ar[rr]^{\hat{h}}_{((A\times\s)\circ\ell,h)}
\ar[d]_\err
\ar @{}[dl]|(0.3)=
\ar@/_1pc/[drr]_h
\ar@{}[drr]|(.5)=
&&(A\times\N)\times B  
\ar[d]_\err
\\
& B
&& B
}
\end{align*}

In the last previous line the left projection from 
$(A\times\N)\times B$ projects on $A\times\N$\,, and in the 
diagram the right projection projects on $B$\,.

So on the left side we get
\begin{align*}
&f\circ(\id_A,\0):A\to A\times\N\to B\\
&=\err\circ k\circ(\id_A,\0)=\err\circ\hat{g}
=\err\circ((\id_A,\0),g)\\
&=g:A\to B.
\end{align*}

The so constructed map $f$ fulfils the Anchor equation.

On the right side we get
\begin{align*}
&f\circ(A\times\s):A\times\N\to A\times\N\to B\\
&=\err\circ k\circ(A\times\s)=\err\circ\hat{h}\circ k
=h\circ k=h\circ(\ell\circ k,
\err\circ k)\\
&=h\circ(\id_{A\times\N},f)
:A\times\N\to(A\times\N)\times B\to B\\
&\text{by the retraction property $(*)$ above}\\
&\text{and the definition of $f:A\times\N\to B$}.
\end{align*}

The so constructed $f$ also fulfils the Step equation.

Thus the construction of $f=[h^\star g]$ 
satisfies the equations Anchor and Step
\textbf{\qed}

\medskip   
We postulate the uniqueness\index{uniqueness} of $f=[h^\star g]$ as an axiom\index{axiom}, 
the last one of theory $\PR$\,.

\subsubsection*{Uniqueness axiom (pr!)}

\inference{(\mr{pr}!)}
{$f':A\times\N\to B$\\ 
&$g:A\to B,\ h:(A\times\N)\times B\to B$ such that\\
&Anchor: $f'\circ(\id_A,\0)=g:A\to B$\\
&Step: $f'\circ(\id_A\times\s)
:A\times\N\to A\times\N\to B$\\
&$=h\circ(\id_{A\times\N},f')
:(A\times\N)\to(A\times\N)\times B\to B$
}
{$f'=[h^\star g]: A\times\N\to B$}

\section{Algebraic operations on $\N$ via (pr)}

\begin{itemize}
\item[$\mul$][$\mul$]
Addition\index{addition}

is already defined above as 
$a+n=\s^\S(a,n):\N\times\N\to\N$\,.
\item[$\mul$][$\mul$]
Definition of multiplication\index{multiplication}
\begin{align*}
& a\mul n=\mul(a,n):\N\times\N\to\N.
\end{align*}

In order to define multiplication by the 
schema $(\mr{pr})$ of primitive recursion, take $\N$ for 
both $A$ and $B$ in the schema:

\inference{(\mr{multiply})}
{$g=g(a)=\0:\N\to\N$\\
&$h=h((a,n),b)=b+a:(\N\times\N)\times\N\to\N$
}
{$\mul(a,n):=f(a,n)=[h^\star g](a,n)$\\
&$:\N\times\N\to\N$ with\\
&Anchor: $\mul(a,\0):=f(a,\0)=g(a)=\0$ and\\
&Step: $\mul(a,\s n)=\mul(a,n)+a$\\
&$:=f(a,\s n)=h((a,n),f(a,n))$
}

Analogously exponentiation\index{exponentiation} $a^n:\N\times\N\to\N$
and factorial\index{factorial} $n!:\N\to\N$ are defined
by this schema.

Of special interest is here the introduction of the

\item[$\mul$]
Predecessor\index{predecessor}

We want to define a map $\p:\N\to\N$ such that $\p(\0)=\0$ 

and $\p(\s n)=n=\id_{\N}(n)$\,.
Since $\p$ is not a binary function, we construct an auxiliary 
function

$\p2=\p2(a,n):\N\times\N\to\N$ from which $\p$ 
will be obtained:

Consider the following instance of schema $(\mr{pr})$\,.

\inference{(\mr{\p2})}
{$g=g(a)=\0:\N\to\N$\\
&$h=h((a,n),b)=n$\\
&$=\err\ell((a,n),b):(\N\times\N)
\times\N\to\N\times\N\to\N$
}
{$\p2=\p2(a,n)=[h^\star g](a,n)$\ with\\
&Anchor: $\p2(a,0)=g(a)=\0$\\
&Step: $\p2(a,\s n)=h((a,n),\p2(a,n))=n$}

The predecessor\index{predecessor} $\p=\p(n):\N\to\N$ now can be 
defined by
\begin{align*}
&\p(\0)=\p2(\0,\0)=\0:\one\to\N\ \text{and}\\
&\p(\s(n))=\p2(\s n,\s n)=n\ \text{i.e.}\\
&\p\,\0=\0:\one\to\N\ \text{and}\\
&\p\,\s=\id_{\N}:\N\to\N
\end{align*}

thus $\p$ is a retraction to $\s$\,.

\item[$\mul$]
Truncated subtraction\index{truncated subtraction}

Truncated or cut-off subtraction\index{cut-off subtraction}
$a\dmin n:\N\times\N\to\N$ is defined by
\begin{align*}
&a\dmin n=\p^\S(a,n):\N\times\N\to\N\ \text{that is}\\
&a\dmin\0=a:\N\to\N\\
&a\dmin\s\,n=\p(a\dmin n):\N\times\N\to\N\\
&=(a\dmin n)\dmin\ttt{1}:\N\times\N\to\N
\end{align*}
\end{itemize}

\section{Boolean operations on $\N$}

\begin{itemize}
\item[$\mul$]
\emph{Boolean truth values}\index{truth value} are defined on 
$\N$ by
\begin{align*}
\false&=\0,\ \true=\ttt{1}:\one\to\N
\end{align*}

Alternatively we also make use of the maps
\begin{align*}
\false_{\N}&=\0=\0\circ\Pi=\false\circ\Pi
:\N\to\one\to\N\\
\true_{\N}&=\ttt{1}=\ttt{1}\circ\Pi=\true\circ\Pi
:\N\to\one\to\N
\end{align*}

\item[$\mul$]
\emph{Negation}\index{negation} is defined on $\N$ by
\begin{align*}
&\neg\alpha=\neg(\alpha)=\ttt{1}\dmin\alpha:\N\to\N
\end{align*}

\item[$\mul$]
\emph{Signum}\index{signum} is defined on $\N$ by
\begin{align*}
&\sg(n)=\neg\neg(n)=\ttt{1}\dmin(\ttt{1}\dmin n):\N\to\N
\end{align*}
 
and hence characterised by
\begin{align*}
\sg(\0)&=\0=\0\circ\Pi:\N\to\one\to\N\\
\sg(\s n)&=\ttt{1}=\ttt{1}\circ\Pi
=\s\circ\0\circ\Pi:\N\to\one\to\N\to\N
\end{align*}

\item[$\mul$]
The logical \emph{and} is defined on $\N$ by
\begin{align*}
&\alpha\land\beta=\sg(\alpha\mul\beta):\N\times\N\to\N
\end{align*}

\item[$\mul$]
The logical \emph{or} is defined on $\N$ by
\begin{align*}
&\alpha\lor\beta=\sg(\alpha+\beta):\N\times\N\to\N
\end{align*}

\item[$\mul$]
The \emph{implication}\index{implication} is defined on $\N$ by
\begin{align*}
&(\alpha\impl\beta)=\neg\alpha\lor\beta:\N\times\N\to\N
\end{align*}
\end{itemize}

\section{Linear order on $\N$}
\begin{itemize}
\item[$\mul$]
\emph{Weak order}\index{weak order} (here a \emph{predicate}) is defined
on $\N$ by
\begin{align*}
&(a\leq b)=(a\impl b):\N\times\N\to\N
\end{align*}

\item[$\mul$]
\emph{Strict order\index{strict order}} is defined on $\N$ by
\begin{align*}
&(a<b)=\sg(b\dmin a):\N\times\N\to\N
\end{align*}

\item[$\mul$]
\emph{Equality\index{equality}} on $\N$ is defined by
\begin{align*}
&(a\doteq b)=(a\leq b)\land(b\leq a):\N\times\N\to\N
\end{align*}

\item[$\mul$]
\paragraph*{Lemma.} Weak order `$\leq$' is reflexive\index{reflexive}
\begin{align*}
&a\leq a
\end{align*}

\item[$\mul$]
\paragraph*{Lemma.} Weak order is antisymmetric\index{antisymmetric}
\begin{align*}
&(a\leq b)\land(b\leq a)\impl(a\doteq b)
\end{align*}

\item[$\mul$]
\paragraph*{Lemma.} Weak order is transitive\index{transitive}
\begin{align*}
&(a\leq b)\land(b\leq c)\impl(a\leq c)
\end{align*}

\item[$\mul$]
\paragraph*{Lemma.} Strict order `$<$' is irreflexive\index{irreflexive}
\begin{align*}
&\neg(a<a)\ \text{and}\\
&  transitive:\\
&(a<b)\land(b<c)\impl(a<c)
\end{align*}

\item[$\mul$] 
\paragraph*{Lemma.} The predicate $(a\doteq b)$ is an
equivalence predicate\index{equivalence predicate} since it is reflexive,
symmetric, and transitive.

\item[$\mul$]
\paragraph*{Lemma.} The law of trichotomy\index{trichotomy} holds 
for $<,\ \doteq,\ >$:
\begin{align*}
&((a<b)\,\dot{\lor}\,(a\doteq b))\,\dot{\lor}\,(b<a)
\end{align*}

where `$\dot{\lor}$' is the \emph{exclusive ``or''},
$\mr{XOR}$\,,
\begin{align*}
&a\,\dot{\lor}\,b=(a\lor b)\land\neg(a\land b)
\end{align*}
\end{itemize}

\paragraph*{Proof} of the Lemmata by the 
truth table method\index{truth table method}\footnote{
  see section 4.4 below on Junctors and truth tables}.

\section{Coproduct construction}

The notion of a sum/coproduct\index{sum/coproduct} with injections 
\begin{align*}
&\xymatrix{
A\ar[r]^(.4)\iota 
&A+B
&B\ar[l]_(.4)\kappa
}
\end{align*}

is \emph{dual} -- revert all arrows --
to the notion of a cartesian -- universal\index{universal} -- product
with projections\index{projection}
\begin{align*}
&\xymatrix{
A&A\times B\ar[l]_(.55)\ell
\ar[r]^(.55)\err
&B
}
\end{align*}

For an object $C$ and maps $f:A\to C,\ g:B\to C$
the couniversal\index{couniversal} \index{couniversal} property of the coproduct is expressed
by the \textsc{diagram}

\begin{align*}
&\xymatrix@+0em{
A
\ar@<0ex>[d]^\iota
\ar[drr]^f
&
\ar@{}[dl]|(.55)=
\\
A+B
\ar@{-->}[rr]^(.45){h}
\ar@{-->}[rr]_(.45){f\choose g}
&&C
\\
B
\ar@<0ex>[u]_\kappa
\ar[urr]_g
&
\ar@{}[ul]|(.55)=
}
\end{align*}

with unique\index{unique}
  
coinduced\index{coinduced} map $h={f\choose g}:A+B\to C$
which is preliminarily to be understood as
\begin{align*}
&h=h(x)=
\begin{cases}
f(a)\ \textbf{if}\ x=\iota(a)\\
\quad \textbf{or else}\\
g(b)\ \textbf{if}\ x=\kappa(b)
\end{cases}: A+B\to C
\end{align*}

\section{Special (equational) coproducts in $\PR$}

\subsubsection*{Theorem}
In $\PR$ the object $\N$\,, together with the 
maps $\0:\one\to\N$
and $\s:\N\to\N$ as injections\index{injection}, is a weak\index{weak}, an \emph{equational\index{equational}}
sum/coproduct: uniqueness of the commutative fill-in
is not required.
\paragraph*{Proof.} We show that $\N$ generates from 
arrows $\one\xto{f}C$
and $\N\xto{g}C$ into an arbitrary object $C$
a coinduced \index{coinduced} map $h={f\choose g}$ into $C$ such that
\begin{align*}
&\xymatrix@+0em{
\one
\ar@<0ex>[d]^\0
\ar[drr]^f
&
\ar@{}[dl]|(.55)=
\\
\N
\ar@{-->}[rr]^(.45){h}
\ar@{-->}[rr]_(.45){f\choose g}
&&C
\\
\N
\ar@<0ex>[u]_\s
\ar[urr]_g
&
\ar@{}[ul]|(.55)=
}
\end{align*} 

Let namely $h=h(m):\N\to C$ be defined 
\emph{arithmetically\index{arithmetically}} by
\begin{align*}
&h=h(m):=\neg(m)\mul f(\0)+\sg(m)\mul g(\p\,m):\N\to C
\end{align*}

where `$\mul$' and `$+$' are defined componentwise\index{componentwise}
for $\PR$ object $C$ which comes as a 
cartesian product of $\N$'s and $\one$'s.

As is seen, when $\iota=\0$ and $\kappa=\s$\,, the 
preliminary definition above generates the desired 
coinduced\index{coinduced} map correctly.

\subsubsection*{Theorem}
For any object $A$\ $\PR$ has the
sum/coproduct\index{sum/coproduct} of object $A$ and object $A\times\N$ in the form
\begin{align*}
&\xymatrix@+1em{
A
\ar[r]^(.4){(\id_A,\0_A)}
&A\times\N
&A\times\N
\ar[l]_(.45){\id_A\times\s}
}
\end{align*}
with sum injections\index{sum injections} 
$\iota=(\id_A,\0_A)=(\id_A,\0\circ\Pi_A):A\to A\times\N$
and $\kappa=\id_A\times\s:A\times\N\to A\times\N$\,.

\paragraph*{Proof.}
We show that $A\times\N$ generates from 
arrows $A\xto{f}C$
and $A\times\N\xto{g}C$ into an arbitrary object $C$
a unique coinduced\index{coinduced} map $h={f\choose g}$ into $C$ such 
that
\begin{align*}
&\xymatrix@+.0em{
A
\ar@<0ex>[d]_{(\id_A,\0_A)}
\ar[drr]^f
&
\ar@{}[dl]|(.55)=
\\
A\times\N
\ar@{-->}[rr]^(.45){h}
\ar@{-->}[rr]_(.45){f\choose g}
&&C
\\
A\times\N
\ar@<0ex>[u]^{\id_A\times\s}
\ar[urr]_g
&
\ar@{}[ul]|(.55)=
}
\end{align*} 

Let namely $h:A\times\N\to C$ be defined 
\emph{arithmetically} by
\begin{align*}
&h=h(a,m):=\neg(m)\mul f(a,\0)+\sg(m)\mul g(a,\p\,m)
:A\times\N\to C
\end{align*}

where `$\mul$' and `$+$' are defined componentwise
for $\PR$ object $C$ which comes as a 
cartesian product of $\N$'s and $\one$'s.

As is seen, when $\iota=(a,\0_A)$ and 
$\kappa=a\times\s$\,,
the definition above generates the desired coinduced\index{coinduced} map 
correctly -- and uniquely since uniqueness \textbf{axiom\index{axiom}} 
$(\mr{pr}!)$
is available

\textbf{\qed}

\chapter{Goodstein Arithmetic $\GA$}

The basic \emph{maps} of the \emph{categorical version} $\GA$ 
of Goodstein\index{Goodstein} Arithmetic are the primitive recursive 
maps zero $\0:\one\to\N$\,,
successor $\s:\N\to\N$\,, addition 
$+:\N\times\N\to\N$\,, predecessor $\p:\N\to\N$\,, 
truncated subtraction $\dmin:\N\times\N\to\N$\,, 
and multiplication 
$\mul\,:\N\times\N\to\N$\,.
His basic Equation and Substitution rules are
covered here by reflexivity, symmetry, transitivity
of equality between maps, and substitution
``axiom\index{axiom}s'' \footnote{
  categorically proved versions of Goodstein's axiom\index{axiom}s}

\inference{\mr{Sb}_1}
{$h:A\to B$\\
&$f(b)=g(b):B\to C$
}
{$f\circ h=g\circ h:A\to B\to C$} 
\inference{\mr{Sb}_2}
{$f,\ g:A\to B$\\
&$h:B\to C$
}
{$h\circ f=h\circ g:A\to B\to C$}

\section{Goodstein uniqueness rules in $\PR$}

For to transfer Goodstein\index{Goodstein}'s
$\mr{U_1}$\,-\,$\mr{U_4}$ we translate $x$ into $n$ and
$F$\,, $S$\,, $P$ into $f$\,, $\s$\,, $\p$\,. Then Goodstein\index{Goodstein}'s
first uniqueness rule becomes

\inference{\mr{U}_1^{\mr{G}}}
{ $f(n) = f(\s n)$ }
{ $f(n) = f(\0)$}

For to transfer the property of \emph{uniqueness} we 
apply Freyd's uniqueness schema $(\mr{FR}!)$ which 
demands formally to  introduce a free parameter $a$ 
into these rules. 
 
\inference{\mr{U}_1}
{$f(a,n)=f(a,\s n):A\times\N\to\N$ }
{$f(a,n)=f\circ(a,n)
:A\times\N\xto{\id}A\times\N\xto{f}\N$\\
&$=f(a,\0)=f\circ(a,\0\circ\Pi_\N)$\\
& $:A\times\N\to A\times\one
\xto{\id_A\times\0}A\times\N\xto{f}\N$
}

\paragraph*{Proof of $\mr{U}_1$\,.}
Consider the \textsc{diagram}
\begin{align*}
&
\xymatrix@+1em{
A\times\N
\ar[r]^(.4){(\id,\0)}
\ar[d]^f
\ar[dr]_(0.4){(f,\0)}
\ar@/_5pc/[ddr]_f
&(A\times\N)\times\N
\ar[r]^{(A\times\N)\times\s}
\ar[d]_(0.5){f\times\N}
\ar[dr]^{f\times s}
\ar@{}[dl]|(0.3)=
\ar@{}[dl]|(0.75)=
\ar@/^2pc/[dd]^(0.8)h
&(A\times\N)\times\N
\ar[d]_{f\times\N}
\ar@{}[dl]|(0.3)=
\ar@{}[dl]|(0.7)=
\ar@/^2pc/[dd]_(0.8)h^{\id_\N^\S\,\circ\,(f\times\N)}
\\
\N
\ar[r]_{(\id_\N,\0)}
\ar[dr]_{\id_\N}
&\N\times\N
\ar[r]_{\N\times\s}
\ar[d]_{\id_\N^\S}^(.3){\ell_{\N,\N}}
\ar@{}[dr]|(0.5)=
\ar@{}[dl]|(0.35)=
&\N\times\N
\ar[d]_{\id_\N^\S}
\\
\ar@{}[ur]|(0.25)=
&\N
\ar[r]^{\id_\N}
&\N
}
\end{align*}

The asserted commutativities of this diagram hold, in 
fact:

$\id_\N^\S=\ell_{\N,\N}:\N\times\N\to\N$ holds by
the basic iteration schema $(\mr{EE})$ applied to the
endomorphism $\id_\N:\N\to\N$\,.
All the equalities symbolised by the inserted `$=$' 
signs are verified straight forward by taking into 
account the abbreviation
$\0=\0_B:=\0\circ\Pi_B:B\to\one\xto{\0}\N$
for objects $B$ of $\PR$\,.

By uniqueness schema $(\mr{FR}!)$ there is a unique
commutative fill-in $h:(A\times\N)\times\N\to\N$ in the 
diagram, namely 

$h=\ell_{\N,\N}\circ(f\times\N):
(A\times\N)\times\N\to\N\times\N\xto{\ell}\N$\,.

So commutativity of the \textsc{diagram} 
gives
\begin{align*}
&f(a,n)=f\circ\id=f:A\times\N\to\N\\
&=\ell_{\N,\N}\circ(f\times\N)
\circ(\id_{A\times\N},\0\circ\Pi_{A\times\N})\\
& :A\times\N\to(A\times\N)\times\N
\to\N\times\N\to\N\\
&=\id_\N\circ f:A\times\N\xto{f}\N\to\N\\
&=\ell_{\N,\N}
\circ(f\times\N)\circ(\id_{A\times\N},0\circ\Pi_\N)\\
& :A\times\N\to(A\times\N)\times\N
\xto{(f\times\N)}\N\times\N\xto{\ell_{\N,\N}}\N\\
&=f\circ(\id_A,\0\circ\Pi_\N):A\times\N
\to A\times\N\xto{f}\N\\
&=f(a,\0):A\times\N\to\N
\end{align*}

\textbf{\qed}
  
\inference{\mr{U}_2}
{$f(a,\s n)=\s f(a,n):A\times\N\to\N$}
{$f(a,n)=f(a,\0)+n:A\times\N\to\N$}\footnote{
In \NAME{Goodstein} 1971, p.\ 114, 
Uniqueness rule $\mr{U_2}$
is misprinted as \inference{\mr{U_2}}{$F(x)=SF(x)$}
{$F(x)=F(0)+x$}
The correct version is \inference{\mr{U_2}}{$F(Sx)=SF(x)$}
{$F(x)=F(0)+x$}
  }
 
\paragraph*{Proof.} Consider the \textsc{diagram}

\begin{align*}
&
\xymatrix{ 
&A\times\N
\ar[rr]^{A\times\s}
\ar[dd]_(.4)h_(.6)f^{\s^\S\circ(f(a,\0),\err)}
&&A\times\N
\ar[dd]_(.4)h_(.6)f^{\s^\S\circ(f(a,\0),\err)}
\\
A
\ar[ur]^{(a,\0_A)}
\ar[dr]_{f(a,\0_A)}
\\
&\N
\ar[rr]^{\s}
&&\N
}
\end{align*}

For $h$ being $f=f(a,n)$ this diagram commutes,
and it also commutes for $h$ being
\begin{align*}
&f(a,\0)+n=\s^\S(f(a,\0_A),n)\\
&=\s^\S\circ(f(a,\0\circ\Pi_\N),\err_{A,\N})\\
&:A\times\N\to\N\times\N\xto{\s^\S}\N.\\
& \text{Thus, by uniqueness schema $(\mr{FR}!)$}\\
&f(a,n)=f=f(a,\0)+n
\end{align*}

\textbf{\qed}
 
\inference{\mr{U}_3}
{$f(a,\s n)=\p f(a,n):A\times\N\to\N$}
{$f(a,n)=f(a,\0)\dmin n:A\times\N\to\N$}\footnote{
In \NAME{Goodstein} 1971, p.\ 114, 
Uniqueness rule $\mr{U_3}$
is misprinted as \inference{\mr{U_3}}{$F(x)=PF(x)$}
{$F(x)=F(0)\dmin x$}
The correct version is 
\inference{\mr{U_3}}{$F(Sx)=PF(x)$}{$F(x)=F(0)\dmin x$}
}

\paragraph*{Proof.} Consider the \textsc{diagram}

\begin{align*}
&
\xymatrix{ 
&A\times\N
\ar[rr]^{A\times\s}
\ar[dd]_(.4)h_(.6)f^{\p^\S\circ(f(a,\0),\err)}
&&A\times\N
\ar[dd]_(.4)h_(.6)f^{\p^\S\circ(f(a,\0),\err)}
\\
A
\ar[ur]^{(a,\0_A)}
\ar[dr]_{f(a,\0_A)}
\\
&\N
\ar[rr]^{\s}
&&\N
}
\end{align*}

For $h$ being $f=f(a,n)$ this diagram commutes,
and it also commutes for $h$ being
\begin{align*}
&f(a,\0)\dmin n=\p^\S(f(a,\0_A),n)\\
&=\p^\S\circ(f(a,\0\circ\Pi_\N),\err_{A,\N})\\
&:A\times\N\to\N\times\N\xto{\p^\S}\N.\\
& \text{Thus, by uniqueness schema $(\mr{FR}!)$}\\
&f(a,n)=f=f(a,\0)\dmin n
\end{align*}

\textbf{\qed}
 
\inference{\mr{U}_4}
{$f(a,\0_A)=g(a,\0_A):A\to\N$\\ 
&$f(a,\s n)=g(a,\s n):A\times\N\to\N$}
{$f(a,n)=g(a,n):A\times\N\to\N$}

\paragraph*{Proof.}
The sum/coproduct of object $A$ and object $A\times\N$ can be given the form
\begin{align*}
&\xymatrix@+1em{
A
\ar[r]^(.4){(\id_A,\0_A)}
&A\times\N
&A\times\N
\ar[l]_(.45){\id_A\times\s}
}
\end{align*}

with sum injections $\iota=(\id_A,\0_A):A\to A\times\N$
and $\kappa=\id_A\times\s:A\times\N\to A\times\N$\,,
cf.\ section 2.8 \emph{special coproducts in $\PR$}.
 
So we have the 
commutative \textsc{diagram}
\begin{align*}
&\xymatrix@+1em{
A
\ar[d]_{(\id_A,\0_A)}
\ar[drr]^{f(a,\0_A)}_{g(a,\0_A)}
\\
A\times\N
\ar@{-->}[rr]^f_g
&&\N
\\
A\times\N
\ar[u]^{\id_A\times\s}
\ar[urr]_{f(a,\s\,n)}^{g(a,\s\,n)}
}
\end{align*}

Thus the two equations in the antecedens of $\mr{U}_4$
entail

$f(a,n)=f=g=g(a,n)$
 
\textbf{\qed}

Categorical Goodstein\index{Goodstein}'s Free-Variables 
Arithmetic $\GA$ and corresponding proofs are thus 
integrated into $\PR$\,.

\section{Deduction of Arithmetic}

We list here -- with 
\emph{passive parameters} $a$ (and others) added 
-- \NAME{Goodstein}'s arithmetical equations
together with his numbering 

1., 1.1, \ldots , 4.4 
referenced as [1.], [1.1], \ldots , [4.4], important of 
them with their 
categorical proofs.

\begin{itemize}
\item[$\mul$]
$
[1.]\quad (a\dmin n)\dmin\ttt{1}
=(a\dmin\ttt{1})\dmin n:\N\times\N\to\N
$

\item[$\mul$] 
$
[1.1]\quad \s a\dmin\s b=a\dmin b:\N\times\N\to\N
$

\item[$\mul$]
$
[1.2]\quad  a\dmin a=\0:\N\to\N
$

\item[$\mul$]
$
[1.3]\quad \0\dmin a=\0:\N\to\N
$

\item[$\mul$] 
$
[1.31]\quad  a\dmin(b+c)=(a\dmin b)\dmin c:
(\N\times\N)\times\N\to\N
$

\item[$\mul$]
$
[1.4]\quad (a+n)\dmin(b+n)
=a\dmin b:(\N\times\N)\times\N\to\N
$

\item[$\mul$] 
$
[2]\quad \0+n = n:\N\to\N
$

\item[$\mul$]
$
[2.1]\quad  a+\s b=\s a+b
$

\item[$\mul$] 
$
[2.2]\quad  a+b=b+a:\N\times\N\to\N
$

\paragraph*{Proof.}

\inference{(\mr{U_4})}
{$a+\0=a=\0+a$ by [2] above\\
&$a+\s b=\s a+b$ by [2.1] above
}
{$a+b=b+a$ \textbf{\qed}}

\item[$\mul$]
$
[2.21] (a\dmin b)\dmin c=(a\dmin c)\dmin b:
(\N\times\N)\times\N\to\N
$

\item[$\mul$]
$
[2.3]\quad (a+n)\dmin n=(a+n)\dmin(\0+n)=a:
\N\times\N\to\N
$ 

\item[$\mul$]
$
[2.4]\quad (a+b)+c=a+(b+c):(\N\times\N)\times\N\to\N
$

\paragraph*{Proof.} For 
$f((a,b),c):=a+(b+c):(\N\times\N)\times\N\to\N$

\inference{(\mr{U_2})}
{$f((a,b),\s c)=a+(b+\s c)=a+\s(b+c)$\\
&$=\s(a+(b+c))=\s f((a,b),c)$
}
{$a+(b+c)=f((a,b),c)=f((a,b),\0)+c$\\
&$=(a+(b+\0))+c=(a+b)+c$ \textbf{\qed}
}

\item[$\mul$]
$ 
[2.5]\quad \0\mul n=\0:\N\to\N
$

\item[$\mul$]
$
[2.6]\quad  a\mul(\ttt{1}\dmin a)=\0
$
\end{itemize}

For the proofs of the properties of multiplication 
we need the
following key theorem, which in 
\NAME{Goodstein} 1957/1964 still comes as an 
axiom\index{axiom}. 
 
\section{Maximum symmetry}

To show the symmetry of the maximum 
algebraically -- that is without making use of the order relation -- is a somewhat intricate endeavour Goodstein\index{Goodstein} submits
himself to, starting with the definition
\begin{align*}
&\max(a,b):=a+(b\dmin a):\N\times\N\to\N
\end{align*}

and proving the

\subsubsection*{Maximum symmetry Theorem}
\begin{align*}
&[3]\quad  a+(b\dmin a)\ =\ b+(a\dmin b)\\
&\quad  \text{that is}\\
&\max(a,b)=\max(b,a):\N\times\N\to\N.
\end{align*}

\paragraph*{Proof} in several steps:
\begin{itemize}
\item
\subsubsection*{Diagonal Reduction Lemma for maximum}

$\max(a,b)=\max(a\dmin\ttt{1},b\dmin\ttt{1})+\sg(a+b)
\ \text{i.e.}$

$[3.1]  a+(b\dmin a)=(a\dmin\ttt{1})+((b\dmin\ttt{1})
\dmin(a\dmin\ttt{1}))+(\ttt{1}\dmin(\ttt{1}\dmin(a+b))$

\paragraph*{Proof} by case distinction:
\begin{align*}
&\max(a,b)\\
&=\bs{if}
\,[(a\doteq\0)\land
(b\doteq\0),(\0\,|\max(a\dmin\ttt{1},b\dmin\ttt{1})
+\sg(a+b))]\\
&= 
\begin{cases}
\0\ \text{if}\ (a\doteq\0)\land(b\doteq\0)\\
\max(a\dmin\ttt{1},b\dmin\ttt{1})+\sg(a+b)\ 
\text{if}\ [\sg(a)
\doteq\ttt{1}]\lor[\sg(b)\doteq\ttt{1}]
\end{cases}\\
&=\max(a\dmin\ttt{1},b\dmin\ttt{1})+\sg(a+b)\\
&\textbf{\qed}
\end{align*}

\item
For proving the symmetry of maximum let
\begin{align*}
&\phi=\phi(n,c):\N\times\N\to\N\ \text{be defined by}\\
&\phi(\0,c)=\0=\0\Pi_{\N}:\N\to\one\to\N\ \text{and}\\
&\phi(\s n,c)=\phi(n,c)+\sg(c\dmin n)\\
&=\sg(c\dmin n)+\phi(n,c)\ \text{by [2.2]}
\end{align*}

\item
For this $\phi$ we \text{show}
\begin{align*}
&\max(a\dmin n,b\dmin n)+\phi(n,a+b)\\ 
&=\max(a\dmin\s n,b\dmin\s n)+\phi(\s n,a+b)&(*)
\end{align*}

\paragraph*{Proof.} Substitution of $(a\dmin n)$
for $a$\,, and $(b\dmin n)$ for $b$ within 
Diagonal Reduction Lemma for maximum gives
\begin{align*}
&\max(a\dmin n,b\dmin n)\\
&=\max((a\dmin n)\dmin\ttt{1},(b\dmin n)\dmin\ttt{1})
+\sg((a\dmin n)+(b\dmin n))\\
&=\max(a\dmin\s n,b\dmin\s n)+\sg((a\dmin n)+(b\dmin n))
\end{align*}

Adding $\phi(n,a+b)$ to both sides of this equation gives
\begin{align*}
&\max(a\dmin n,b\dmin n)+\phi(n,a+b)\\
&=\max(a\dmin\s n,b\dmin\s n)+\sg((a\dmin n)+(b\dmin n))\\
& +\phi(n,a+b)\\
&=\max(a\dmin\s n,b\dmin\s n)+\phi(\s n,a+b)
\ \text{by def.\ of $\phi$}
\end{align*} 

\item
We also \text{show} for this $\phi$
\begin{align*}
&\max(b\dmin n,a\dmin n)+\phi(n,b+a)\\ 
& = \max(b\dmin\s n,a\dmin\s n)+\phi(\s n,b+a)&(**)
\end{align*}

\paragraph*{Proof.}
Exchange $a$ and $b$ in $(*)$ above.

\item
In order finally to \text{prove} the 
symmetry of the binary maximum we 

apply $\mr{U}_1$ and the Diagonal Reduction Lemma.
\begin{align*}
&\max(a\dmin n,b\dmin n)+\phi(n,a+b)\\
&=\max(a\dmin\0,b\dmin\0)+\phi(\0,a+b)
\ \text{by $(**)$ and $\mr{U}_1$}\\
&=\max(a,b)+\0\ \text{by Red.\ Lemma}\\
&=\max(a,b)\\
& \text{as well as}\\
&\max(b\dmin n,a\dmin n)+\phi(n,b+a)\\
&=\max(b\dmin\0,a\dmin\0)+\phi(\0,b+a)\\
&=\max(b,a)+\0\\
&=\max(b,a)
\end{align*}

$(*)$ and $(**)$ put together:
\begin{align*}
&\max(a,b)=\max(a\dmin n,b\dmin n)+\phi(n,a+b)\\
& \text{and}\\
&\max(b,a)=\max(b\dmin n,a\dmin n)+\phi(n,b+a)\\
&=\max(b\dmin n,a\dmin n)+\phi(n,a+b)\ \text{by [2.2]}
\end{align*}

and \text{hence}, by substituting free variable $n$ 
by $b$ in the 

latter two equations:
\begin{align*}
&\max(a,b)=\max(a\dmin b,b\dmin b)+\phi(b,a+b)\\
&=\max(a\dmin b,\0)+\phi(b,a+b)\\
&=(a\dmin b)+\phi(b,a+b)&(\bs{\circ})
\end{align*}

\textbf{and}
\begin{align*}
&\max(b,a)=\max(b\dmin b,a\dmin b)+\phi(b,a+b)\\
&=\max(b\dmin b,a\dmin b)+\phi(b,a+b)\\
&=\max(\0,a\dmin b)+\phi(b,a+b)\\
&=(a\dmin b)+\phi(b,a+b)&(\bs{\circ}\bs{\circ})
\end{align*}

\text{hence} finally, by $(\bs{\circ})$ and 
$(\bs{\circ}\bs{\circ})$\,, 

and symmetry and transitivity of `$=$'
\begin{align*}
&\max(a,b)=\max(b,a):\N\times\N\to\N
\end{align*}
\end{itemize}
\textbf{\qed} symmetry of the binary maximum.

\medskip\noindent  
From his axiom\index{axiom}s (here theorems) Goodstein\index{Goodstein} 
proves a schema named 
$\mr{V}_4$ in 
correspondence to $\mr{U}_4$ which is very helpful for 
various purposes. 
   
We \text{show} $\mr{V}_4$ for $\GA$ 
(and hence for $\PR$):

\inference{\mr{V}_4}
{$f, g, h: A\times\N\to\N$\\
&$f(a,\0)=g(a,\0):A\to\N$ \\
&$f(a,\s\,n)=f(a,n)+h(a,n):A \times\N\to\N$\\ 
&$g(a,\s n)=g(a,n)+h(a,n):A\times\N\to\N$  
}
{$f(a,n)=g(a,n)$}

\begin{minipage}{\textwidth}
\paragraph*{Proof.}
Rule $\mr{V}_4$ can be \text{deduced} by applying 
$\mr{U}_1$ to the distance map 
\begin{align*}
&|f(a,n)-g(a,n)|
=(f(a,n)\dmin g(a,n))+(g(a,n)\dmin f(a,n))\\
&:A\times\N\to\N
\end{align*}
as follows:
\begin{align*} 
&(f(a,\s n)\dmin g(a,\s n))+(g(a,\s n)\dmin f(a,\s n))\\ 
&=(f(a,n)+h(a,n))\dmin(g(a,n)+h(a,n))\\ 
& +(g(a,n)+h(a,n))\dmin(f(a,n)+h(a,n))\\
& \text{by antecedent of $\mr{V}_4$}\\ 
&=(f(a,n)\dmin g(a,n))+(g(a,n)\dmin f(a,n))\\
&\text{by simplification [1.4]}
\end{align*}
  
So we have by $\mr{U}_1$
\begin{align*}
&(f(a,n)\dmin g(a,n))+(g(a,n)\dmin f(a,n))\\
&=(f(a,\0)\dmin g(a,\0))+(g(a,\0)\dmin f(a,\0))\\
&=\0+\0=\0\ \text{by antecedent 
of $\mr{V}_4$\,, hence}\\
&f(a,n)\dmin g(a,n)=\0=g(a,n)\dmin f(a,n):A\times\N \to\N
\end{align*}
    
and \text{hence} by \text{symmetry of maximum} 
\begin{align*} 
f(a,n)&=f(a,n)+(g(a,n)\dmin f(a,n))\\ 
&=\max(f(a,n),g(a,n))\\ 
&=\max(g(a,n),f(a,n))\\
&=g(a,n)+(f(a,n)\dmin g(a,n))\\ 
&=g(a,n)
\end{align*} 

\textbf{\qed} for $\mr{V}_4$
\end{minipage}

\subsubsection*{Distributivity and associativity of multiplication}

\begin{itemize}
\item[$\mul$]
First is shown

$
[4.]\quad  f(a,n):=\s a\mul n=g(a,n):=a\mul n+n
$

using $\mr{V}_4$\,.

\item[$\mul$]
$
[4.2]\quad  a\mul(m+n)=a\mul m+a\mul n
$

\item[$\mul$]
$
[\mr{dis\ r}]\quad  (m+n)\mul a=m\mul a+n\mul a
$

\paragraph*{Proof.}
For $n=\0$ we have by definition of $+$ and \,$\mul$

$
(m+\0)\mul a=m\mul a=m\mul a+\0=m\mul a+\0\mul a.
$

From $(m+n)\mul a=m\mul a+n\mul a$ we conclude to

$   
(m+\s n)\mul a=m\mul a+\s n\mul a
$ 
   
 -- with passive parameter $(m,a)$ -- by $(\mr{V}_4)$
with choice of $f$\,, $g$\,, and $h$ as follows.

We have
\begin{align*} 
f((m,a),n)&:=(m+n)\mul a\\ 
g((m,a),n)&:=m\mul a+n\mul a \mbox{and} \\ 
h((m,a),n)&:=a 
\end{align*}

and
\begin{align*} 
f((m,a),\s n)&=(m+\s n)\mul a=(m+(n+\ttt{1}))\mul a\\ 
&=(m+\s n)\mul a=\s(m+n)\mul a\\ 
&=(m+n)\mul a+a\ \text{by $[4.]$}\\ 
&=f((m,a),n)+h((m,a),n)\\ 
g((m,a),\s n)&=m\mul a+\s n\mul a\\ 
&=m\mul a+n\mul a+a\ \text{by $[4.]$}\\ 
&=g((m,a),n)+h((m,a),n)
\end{align*}

From this $\mr{V}_4$ gives 
\begin{align*}
&f((m,a),n)=g((m,a),n)\ \text{i.e.}\\ 
&(m+n)\mul a=m\mul a +n\mul a\\
&\textbf{\qed}
\end{align*}

\item[$\mul$]
Associativity

$
[4.3]\quad  a\mul(m \mul n)=(a\mul m)\mul n
$

\paragraph*{Proof.}
We prove the law by application of rule $(\mr{V}_4)$ 
with recursion parameter $n$ and passive parameter 
$(a,m)$ to 
\begin{align*} 
f((a,m),n)&:=a\mul(m\mul n)\\ 
g((a,m),n)&:=(a\mul m)\mul n \mbox{and}\\ 
h((a,m),n)&:=a\mul m 
\end{align*}

For $n=\0$ we have $a\mul(m\mul n)=a\mul\0=\0$ 

and on the other hand $(a\mul m)\mul n=(a\mul m)\mul\0=\0$

For $(\mr{V}_4)$ -- step we take
\begin{align*} 
f((a,m),\s n)
&=a\mul(m\mul\s n)\\
&=a\mul(m\mul n+m)\\
&=a\mul(m\mul n)+a\mul m\ \text{by distributivity}\\ 
&=f((a,m),n)+h((a,m),n)\\ 
g((a,m),\s n)
&=(a\mul m)\mul\s n\\
&=(a\mul m)\mul(n+1)\\ 
&=(a\mul m)\mul n+(a\mul m)\mul1
\ \text{by distributivity}\\ 
&=(a\mul m)\mul n+(a\mul m)\\ 
&=g((a,m),n)+h((a,m),n)
\end{align*}

By $(\mr{V}_4)$ we get 
\begin{align*}
&f((a,m),n)=g((a,m),n)\ \text{i.e.}\\ 
&a \mul (m \mul n) = (a \mul m) \mul n
\end{align*}

\textbf{\qed}

\item[$\mul$]
\subsubsection*{Commutativity}
$
[4.1]\quad  a\mul n=n\mul a:\N\times\N\to\N
$

\paragraph*{Proof}
by $(\mr{pr}!)$\footnote{
  here we don't follow Goodstein\index{Goodstein}'s proof}.
\begin{align*}  
&\text{For}\ \,f(a,n)=a\mul n\ \,\text{and}\ \, 
f'(a,n)=n\mul a:\N\times\N\to\N\\
&\text{and}\ \,h=h((a,n),b)=b+a:(\N\times\N)\times\N\to\N
\ \,\text{we have}\\
&f(a,\0)=a\mul\0=\0=\0\mul a=f'(a,\0)
\ \,\text{as well as}\\
&f(a,\s\,n)=a\mul\s\,n\\
&=a\mul(n+\ttt{1})=a\mul n+a\mul\ttt{1}
\ \,\text{by distributivity}\\
&=a\mul n+a=h((a,n),f(a,n))\ \,\text{and}\\
&f'(a,\s n)=\s n\mul a=(n+\ttt{1})\mul a\\
&=n\mul a+\ttt{1}\mul a\ \,\text{by $[\mr{dis\ r}]$ above}\\
&=n\mul a+a=h((a,n),f'(a,n)).
\end{align*}

So $f(a,n)=f'(a,n)$\ \text{i.e.}
$a\mul n=n\mul a:\N\times\N\to\N$ 

follows by the uniqueness schema $(\mr{pr}!)$ of primitive
recursion

\textbf{\qed}

\item[$\mul$]
$
[4.4]\quad  a\mul(m \dmin n)=a\mul m\dmin a\mul n
$

\paragraph*{Proof} by 
$(\mr{V}_4)$ -- with passive parameter $(a,m)$\,. 
\end{itemize}

As remarked in the Introduction, Goodstein\index{Goodstein} also constructs
formal logic out of his logic-free recursive Arithmetic
including modus ponens, deriving $Q$ from $P$ and
$P\to Q$\,. This, he shows, corresponds 
(``is associated'') to getting $q=0$ from
$p=0$ and $(1\dmin p)q=0$ which he proves to be valid in
his recursive Arithmetic. 

The schema of induction\index{induction}, to infer $P(x)$ from $P(0)$
and $P(x)\to P(Sx)$\,, Goodstein\index{Goodstein} gets from
modus ponens by a suitable recursive function $q(x)$\,,
cf.\ \NAME{Goodstein} 1971, pp.\ 120-121.

\chapter{Subobjects by predicates. Theory $\bfS$}

If $\PR$ is supplemented by 
predicate comprehension schemata, a wider 
realm of categorical subobjects 
arises which we call
Skolem\index{Skolem} Theory $\bfS$ here.

\section{Definitional extension of $\PR$ to $\bfS$}

First we introduce the predicate-into-subobject 
\emph{comprehension} schema

\inference {(\bfS_{\Obj})}
{$\chi:A\to\N$ $\PR$ predicate}
{$\set{A:\chi}$ object of theory $\bfS$}
  
For $\set{A:\chi}$ we may also write
$\set{a\in A:\chi(a)}$
where $a$ is a \emph{bound variable}.
The $\in$ sign is used here informally without
commitment to axiomatic set theory.

The $\bfS$-\emph{maps} from subobject
$\set{A:\chi}$ of $A$ to subobject
$\set{B:\ph}$ 
of $B$ are those $\PR$ maps $f:A\to B$ which restrict 
in $\bfS$ to 
$\set{A:\chi}$ as domain and to $\set{B:\ph}$ as codomain.

We introduce the following \text{definition schema} for 
map-extension 

\inference { (\Ext_{\mathbf{Map}}) }
{$\set{A:\chi},\ \set{B:\ph}$ $\bfS$-objects,\\
&$f:A\to B$ a $\PR$-map s.t.\\
&$\PR\derives(\chi(a)\impl\ph f(a))$
}
{$f$ is an $\bfS$-map $f:\set{A:\chi}\to\set{B:\ph}$}
 
\medskip
In case $A=B$ and $\chi', \chi'':A\to\N$ predicates, if 
\begin{align*}
& \PR\derives[\chi'(a)\impl\chi''(a)]:A\to\N,
\end{align*}

then the $\bfS$-map
\begin{align*}
&\id_A:\set{A:\chi'}\to\set{A:\chi''}\\
& \text{is called an inclusion, and written}\\ 
&\subseteq\,:\set{A:\chi'}\to\set{A:\chi''}\ \text{or}\ 
\set{A:\chi'}\subseteq\set{A:\chi''}.
\end{align*}

\paragraph*{Observation:} 
For predicates $\chi, \ph: A \to \N$ such that
$\PR \derives\chi = \ph: A \to \N$ holds
\begin{align*}
& \set{A:\chi} \subseteq \set{A:\ph}\ \text{and}
      \ \,\set{A:\ph} \subseteq \set{A:\chi}.
\end{align*}

Rather than implying the equality of $\set{A:\chi}$ 
and $\set{A:\ph}$
this amounts only to an isomorphism
\begin{align*}
&\id_A: \set{A:\chi} \xto{\iso} \set{A:\ph}\,.
\end{align*}

This endows the construction of $\bfS$ with an inherent 
and general property
wherein identity is not defined by the property of 
containing the same elements
but only by equivalence to an identity map.

In particular an \emph{inclusion}
\begin{align*}
&\subseteq\,: \set{A:\chi'}\to\set{A:\chi''}
\end{align*}

is only an inclusion up to isomorphism.

The equality predicate
$(m\doteq n):\N\times\N\to\N$ on $\N$ and then on a 
general $\PR$-object $A$\,, is carried over to $\bfS$ as 
follows.
\begin{align*}
&(m\leq n)=\neg(m\dmin n):\N\times\N\to\N\\
&(m\doteq n)
=\neg((m\dmin n)+(n\dmin m))\\
&=(m\leq n)\land(n\leq m):\N\times\N\to\N\\
& \text{
and further recursively, componentwise}\\
&((a,b)\doteq(a',b'))=(a\doteq a')\land(b\doteq b')\\
&:(A\times B)\times(A\times B)
\xto{(\ell\times\ell)\times(\err\times\err)}(A\times A)
\times(B\times B)\\
&\xto{\doteq\times\doteq}(\N\times\N)
\xto{\land}\N.
\end{align*}

We introduce equality `$=$' of \emph{maps} in $\bfS$ by the 
schema

\inference {(\bfS_{=})}
{$f, g:\set{A:\chi}\to\set{B:\ph}$ $\bfS$-maps\\
&$\PR\derives\chi(a)\impl[f(a)\doteq g(a)]$  
}
{$\bfS\derives f=g:\set{A:\chi}\to\set{B:\ph}$}

From equality of maps in $\PR$ reflexivity, symmetry, 
and transitivity is inherited.

\section{Cartesian category structure of $\bfS$}

(For this cf.\ \NAME{Eilenberg, Elgot} 1970 I.3 and I.4.)

The \emph{identity map} for an 
$\bfS$ object $\set{A:\chi}$ 
is simply
\begin{align*}
&\id_{\set{A:\chi}}=\id_A:\set{A:\chi}\to\set{A:\chi}
\end{align*}

The \emph{composition} of the $\bfS$-maps 
\begin{align*}
&f:\set{A:\chi}\to\set{B:\ph}\ \,\text{and}\ \,
g:\set{B:\ph}\to\set{C:\psi}
\end{align*} 

is inherited from $\PR$ as well:
\begin{align*}
&\derives(\chi(a)\impl\ph(f(a)))\ \text{and}\ \,
\derives(\ph(b)\impl\psi(g(b)))\ \text{hence}\\ 
&\derives(\chi(a)\impl\ph(f(a))\impl\psi(g(f(a))))\,.
\end{align*} 

The map 
\begin{align*}
&g\circ f:\set{A:\chi}\to\set{B:\ph}\to\set{C:\psi}
\end{align*}

becomes a composition map in $\bfS$ 
by transitivity of `$\impl$'. The neutrality of 
identity maps, and the associativity of composition 
are inferred too.

The cartesian product is defined in $\bfS$ by
\begin{align*}
&\set{A:\chi}\times\set{B:\ph}
:=\set{A\times B:\chi\land\ph}\,.
\end{align*}

This becomes more transparent when written
\begin{align*}
&\set{a\in A:\chi(a)}\times\set{b\in B:\ph(b)}\\
&=\set{(a,b)\in A\times B:\chi(a)\land\ph(b)}\,,
\end{align*}
  
\textsc{diagram}
\begin{align*}
&\xymatrix{\set{A:\chi}
&\set{A:\chi}\times\set{B:\ph}
\ar[l]_(.6)\ell
\ar[r]^(.6)\err
&\set{B:\ph}
}
\end{align*}

In this diagram the $\bfS$-product is 
constructed by domain/codomains restriction of the 
corresponding $\PR$ 
cartesian product projections, \textsc{diagram}:
\begin{align*}
\xymatrix{
&A
&A\times B
\ar[l]_\ell
\ar[r]^\err
&B
}
\end{align*}

\subsubsection*{Theorem}
The universal\index{property} property of the cartesian product in $\bfS$ is 
given by forming the unique induced map\index{induced map} which in turn is 
defined by domain/co\-domain restriction of the induced map\index{induced map} 
in $\PR$\,.

\paragraph*{Proof.}
Consider the induced-into-product 
\textsc{diagram}
\begin{align*}
&\xymatrix@+0em{
&\set{A:\chi}
\\
\set{C:\psi}
\ar@/^1pc/[ru]^{f}
\ar@{}[ru]|{\overset{\phantom{M}}{=}} 
\ar@{-->}[r]^(0.4){(f,g)}
\ar@{}[rd]|{\underset{\phantom{M}}{=}} 
\ar@/_1pc/[rd]_{g}
&\set{A\times B:\chi\land\ph}
\ar [u]_{\ell}
\ar [d]^{\err}
\\
&\set{B:\ph}
}
\end{align*}

Here $f$ and $g$ are assumed to be $\bfS$-maps, and the 
fill-in $(f,g)$ is to be shown to be an $\bfS$-map too.

In $\PR$ we have
\begin{align*}
&\derives\psi(c)\impl\chi\,f(c)\ \text{and}\ 
\derives\psi(c)\impl\ph\,g(c)\,,\,\text{hence}\\
&\derives(\psi(c)\impl\chi\,f(c))\land(\psi(c)
\impl\ph\,g(c))\,,\\
& \text{and therefore}\\
&\derives\psi(c)\impl\chi\,f(c)\land\ph\,g(c)\,.
\end{align*}

The latter derivation reveals the unique $\PR$-fill-in
\begin{align*}
&(f,g):C\to A\times B\\
& \text{as the (unique) $\bfS$-fill-in}\\
&(f,g):\set{C:\psi}\to\set{A\times B:\chi\land\ph}
\end{align*}

\textbf{\qed}

\subsubsection*{Definition} of $\bfS$-predicates.

We define an $\bfS$-predicate to be an $\bfS$-map
\begin{align*} 
&\psi=\psi(a):\set{A:\chi}\to\set{\N:\true_{\N}}
\ \text{with}\ \psi=\sg\circ\psi\,.
\end{align*}

For $\chi=\chi(a):A\to\N$ a $\PR$ predicate, the $\PR$ 
predicative equality
\begin{align*}
&(a\doteq a'):A\times A\to\N
\end{align*}

restricts to a (reflexive, symmetric, transitive)
$\bfS$-equality
\begin{align*}
&(a\doteq a'):\set{A:\chi}\times\set{A:\chi}\to\N
\end{align*}

which is preserved by any 
$\bfS$-map $f:\set{A:\chi}\to\set{B:\ph}:$
\begin{align*}
&\bfS\derives(a\doteq a')\impl(f(a)\doteq f(a'))
\end{align*}

\subsubsection*{Enlarging $\PR$ to $\bfS$}

$\PR$ is embedded into $\bfS$ as a cartesian 
categorical
subsystem by the Cartesian Inclusion Functor
\begin{align*}
&\mbf{I}:\PR\overset{\sqsubset}{\lto}\bfS\ 
\text{defined by}\\
&[A\xto{f}B]
\mapsto[\set{A:\true_A}\xto{f}\set{B:\true_B}]
\end{align*}

Here $\true_A=\true\circ\Pi_A:
A\xto{\Pi}\one\xto{\true}\N$\,,
 $\true_B:B\to\N$ analogously.

\section{The categorical truth object}
For the categorical \emph{truth object}
\begin{align*}
&\two\defeq\set{\alpha\in\N:\alpha\leq\ttt{1}}\\
&=\set{\alpha\in\N:(\alpha\doteq\0)\dot{\lor}
(\alpha\doteq\ttt{1})}\\
&=\set{\alpha\in\N:(\alpha\doteq\0)+
(\alpha\doteq\ttt{1})}
\end{align*}

holds in $\bfS$ a restricted weak sum/coproduct property,
as we show here:
(co)universality\index{couniversality} \textsc{diagram} is
\begin{align*}
&\xymatrix@+0em{
\one
\ar@<0ex>[d]^\iota
\ar[drr]^f
&
\ar@{}[dl]|(.55)=
\\
\two
\ar@{-->}[rr]^(.45){h}
\ar@{-->}[rr]_(.45){f\choose g}
&&\set{C:\gamma}
\\
\one
\ar@<0ex>[u]_\kappa
\ar[urr]_g
&
\ar@{}[ul]|(.55)=
}
\end{align*}

with $\iota=\0$\,,\,$\kappa=\s\,\0$\,, and with 
componentwise $\dot{\lor}=+$\ 
as well as\ $\land=\,\mul$\,.

\paragraph*{Proof.} 
Because of the closure property
of $\set{C:\gamma}$ we get
\begin{align*} 
&h=h(\alpha)
=((\alpha\doteq\0)\mul f(\0)+(\alpha\doteq\ttt{1})\mul g(\0)):
\two\to\set{C:\gamma}
\end{align*}

as $\bfS$ commutative fill-in into the
\textsc{diagram}.

This shows $\two$ to be the restricted weak coproduct 
\begin{align*}
&\two\iso\one+\one
\end{align*}

with respect to objects $\set{C:\gamma}$ closed under
componentwise
`$\mul$' and `$+$', with injections
\begin{align*}
&\iota=\0:\one\to\two\\ 
&\text{and}\ \kappa=1=\s\,\0:\one\to\two
\end{align*}

\textbf{\qed}

\subsubsection*{Corollary}
$
\xymatrix{\one\ar[r]^\iota_\0&\two&\one\ar[l]_\kappa^1}
$ 
has the
couniversal\index{couniversal} property of a weak 
coproduct\index{coproduct} with respect to any
(finite) cartesian product of $\one$'s, 
$\N$'s -- and $\two$'s.

In the sequel we rely on this Corollary.

\paragraph*{Remark.}For $\bfS$-maps $\one\xto{f}\set{C:\gamma}$\,, $\one\xto{g}\set{C:\gamma}$ the coinduced\index{coinduced} map
$h:\two\to\set{C:\gamma}$ fulfils
\begin{align*}
&h=h(\alpha)=
\begin{cases}
f=f(\0)\ \textbf{if}\ \alpha\doteq\0\\
\quad \textbf{or else}\\
\,g=g(\0)\ \textbf{if}\ \alpha\doteq\ttt{1}
\end{cases}
\end{align*}

\section{Junctors and truth tables}

\subsubsection*{Unary junctors}
A unary \emph{junctor} $\alpha:\two\to\two$ in general 
is given by

its components $\alpha_0, \alpha_1$ of the 
induced map\index{induced map}
\begin{align*}
&\xymatrix@+0em{
\one
\ar[d]_\0
\ar[drr]^{\alpha_0}
\ar@{}[dr]|(.5)=
\\
\two
\ar@{-->}[rr]^(.4){\alpha_0\choose\alpha_1}
&&\two
\\
\one
\ar[u]^{\ttt{1}}
\ar[urr]_{\alpha_1}
\ar@{}[ur]|(.6)=
}      
\end{align*}

\subsubsection*{Unary truth tables}
\emph{Negation} 
$\neg={\ttt{1}\choose\0}:\two\iso\one+\one\to\two$
corresponds to its truth table among

\begin{tabular}[]{c|c}
$\alpha$ &\\
\hline
$\0$ & $\alpha_0$\\
$\ttt{1}$ & $\alpha_1$
\end{tabular} 
\begin{tabular}[]{c|c}
$\sg=\id$ &\\
\hline
$\0$ & $\0$\\
$\ttt{1}$ & $\ttt{1}$
\end{tabular} 
\begin{tabular}[]{c|c}
$\neg$ &\\
\hline
$\0$ & $\ttt{1}$\\
$\ttt{1}$ & $\0$
\end{tabular}

\subsubsection*{Binary junctors}
By the universal\index{universal} properties of the Cartesian product
in $\bfS$ we have 
\begin{align*}
&(\two\times\two)=((\one+\one)\times(\one+\one))\\
&\iso((\one\times\one)+(\one\times\one))+((\one\times\one)+
(\one\times\one))\\
&\iso((\one+\one)+(\one+\one))\\
&=(\two+\two)
\end{align*}

Hence \emph{binary junctors}
$\alpha:(\two\times\two)\to\two$ are given by
their components 
$\alpha_0={\alpha_{00}\choose\alpha_{01}}:\two\to\two$ and
$\alpha_1={\alpha_{10}\choose\alpha_{11}}:\two\to\two$ 
of the coinduced\index{coinduced} map 
$\alpha:\two\times\two\iso\two+\two\to\two$ in the
\textsc{diagram}
\begin{align*}
&\xymatrix@-.7em{
&&\one
\ar@{}[dddl]|(.35)=
\ar[dl]_{\ttt{1}}
\ar[dddr]^{\alpha_{01}}
\\
\ar@{}[ddr]|(.45)=
&\two
\ar[dd]^{(\0,\id)}
\ar[drdr]^{\alpha_0}
\ar@{}[ddr]|(.65)=
\\ 
\one
\ar[ur]^\0
\ar[drrr]_(.2){\alpha_{00}}
\\
&\two\times\two
\ar@{-->}[rr]^(.4)\alpha
&&\two
\\
\one
\ar[dr]_\0
\ar[urrr]^(.2){\alpha_{10}}
\\
\ar@{}[uur]|(.45)=
&\two
\ar[uu]_{(\ttt{1},\id)}
\ar[urur]_{\alpha_1}
\ar@{}[uur]|(.65)=
\\
&&\one
\ar[ul]^{\ttt{1}}
\ar[uuur]_{\alpha_{11}}
\ar@{}[uuul]|(.35)=
}      
\end{align*}
 
\subsubsection*{Truth table}
\begin{tabular}[]{c|c c}
$\alpha$ & $\0$ & $\ttt{1}$\\
\hline
$\0$ & $\alpha_{00}$ & $\alpha_{01}$\\
$\ttt{1}$ & $\alpha_{10}$ & $\alpha_{11}$\\
\end{tabular} 
 
\paragraph*{Cases.} $\land, \lor, \impl$\,, 
\emph{Scheffer stroke} $\bs{|}\ (=\mr{NAND})$\,.

\section{R\'esum\'e}

Taken together we rearive at the \textsc{diagram}
\begin{align*}
&\xymatrix @+0em{
\one
\ar@<+0.5ex>[r]^\bot
\ar @<-0.5ex> [r]_\top 
&\two 
\ar[dr]^(.5)\neg
&&\two\times\two
\ar@<-1.0ex>@/_1pc/[ll]_{\bs{|}}
\ar@<-0.5ex>[ll]_\land
\ar@<+0.5ex>[ll]^\lor
\ar@<-0.5ex>[dl]_(.5)\impl
\ar@<+0.5ex>[dl]^\biimpl
\\
&&\two
}
\end{align*}

with
\begin{itemize}
\item[$\mul$]
truth value false $\bot:\one\to\two$\,,

\item[$\mul$]
truth value true $\top:\one\to\two$\,, \,and

\item[$\mul$]
\emph{Scheffer stroke} $(\mr{NAND})$

$\bs{|}=\alpha\bs{|}\beta=\neg(\alpha\land\beta):
\two\times\two\to\two\to\two$

where the latter generates the boolean 
functions
\begin{align*}
&\neg\alpha=\alpha\bs{|}\alpha\\
&\alpha\lor\beta
=(\alpha\bs{|}\beta)\bs{|}(\alpha\bs{|}\beta)\\
&\alpha\land\beta=\neg(\neg\alpha\lor\neg\beta)\\
&(\alpha\impl\beta)=\neg\alpha\lor\beta\\
&(\alpha\biimpl\beta)
=(\alpha\impl\beta)\land(\beta\impl\alpha)
\end{align*}
\end{itemize}

\subsubsection*{Boolean algebra theorem}
The \emph{truth algebra} 
$\two=\bfan{\two, \bot, \top, \neg, \land, \lor, \impl}$
satisfies the known boolean equations.
\paragraph*{Proof} for unary operators by truth-table 
verification using the two coproduct injections 
$\0:\one\to\two,\ \,\ttt{1}:\one\to\two$\,, and
for binary operators using the four coproduct injections
\begin{align*}
&(\0,\0):\one\to\two\times\two\,,
\ (0,\ttt{1}):\one\to\two\times\two\\
&(\ttt{1},\0):\one\to\two\times\two\,,
\ (\ttt{1},\ttt{1}):\one\to\two\times\two
\end{align*}

by boolean truth-table verification as well\ \,\textbf{\qed}

\subsubsection*{Lemma} 
For $A=\set{A:\alpha}$ in $\bfS$ closed under 
`$\dot{\lor}$' and `$\mul$'\ \text{we have}
$A\times\two=A\times(\one+\one)\iso A+A$\ 
\text{distributively}\,,
closed under componentwise `$\dot{\lor}$' and `$\mul$'\,,
with coproduct injections
\begin{align*}
&\iota=(\id_A,\false_A):A\to A\times\two\ \text{and}\\ 
&\kappa=(\id_A,\true_A):A\to A\times\two
\end{align*}

\textbf{\qed}

\subsubsection*{Case distinction}
For given $g,h:\set{A:\chi}\to\set{C:\gamma}$\,, the 
latter
object closed under componentwise `$\dot{\lor}$' and 
`$\mul$'\,, we introduce a map
\begin{align*}
&f=f(a)=\bs{if}\,[\chi,(g|h)](a)\\
&=
\begin{cases}
g(a)\ \textbf{if}\ \chi(a)\\
\quad \textbf{or else}\\
h(a)\ \textbf{if}\ \neg\,\chi(a)
\end{cases}
:\set{A:\chi}\to\set{C:\gamma}
\end{align*}

This $f$ is characterised by the properties
\begin{align*}
&\chi(a)\impl f(a)\doteq g(a)\ \,\text{and}\\         
&\bs{\neg}\chi(a)\impl f(a)\doteq h(a)       
\end{align*}

\chapter{Order and equality on the NNO}

We begin with the weak and the strict order
predicates, and proceed with the equality predicate.
Their connection with each other and
with the algebraic structure of $\N$ leads in particular
to \emph{equality definability}, a
consequence of Goodstein\index{Goodstein}'s symmetry of the algebraically defined maximum
operation. Finally we state in categorical terms the 
strict resp.\ weak monotony of the binary algebraic 
operations on $\N$

\paragraph*{Definition} 
\begin{itemize}
\item[$\mul$]
A map $\chi=\chi(a):A\to\N$ is 

-- as defined above -- 

a \emph{predicate} on $A$ if
\begin{align*}
&\chi=\sg\,\chi:A\to\N\to\N\ \,\text{where}\\
&\sg(\0)=\0\ \text{and}\ \sg(\s\,n)=1=\s\,\0\,,\ n\
 \text{free in}\ \N\,.
\end{align*}

\item[$\mul$]
Negation is defined by
\begin{align*}
&\neg\,a=1\dmin a:\N\to\N\,.
\end{align*}

\item[$\mul$]
We have $\sg=\sg\,\sg:\N\to\N\to\N$ and 

$\neg=\neg\,\sg=\sg\,\neg:\N\to\N\to\N$\,,

and therefore in particular for any map $f:A\to B$ 
\begin{align*}
&\sg\,f:A\to\N\to\N\ \text{and}\ \neg\,f:A\to\N\to\N
\end{align*}

are predicates.

\item[$\mul$]
Weak order is defined by
\begin{align*}
&(a\leq b)=\neg(a\dmin b)=\N\times\N\to\N\,.
\end{align*}

\item[$\mul$]
Strict order then is defined by
\begin{align*}
&(a<b)=(\s a\leq b)=\sg(b\dmin a):\N\times\N\to\N\,.
\end{align*}

\item[$\mul$]
The \emph{distance} of natural numbers $a,b$ is defined 

by the map
\begin{align*}
&\dist=\dist(a,b)=\abs{a-b}=(a\dmin b)+(b\dmin a):
\N\times\N\to\N\,.
\end{align*}

\item[$\mul$]
The \emph{equality predicate} 
$(a\doteq b):\N\times\N\to\N$ is defined by
\begin{align*}
&(a\doteq b)=(a\leq b)\land(b\leq a)\\
&=\neg(a\dmin b)\land\neg(b\dmin a)\\
&=\neg\,\sg(a\dmin b)\land\neg\,\sg(b\dmin a)\\
&=\neg(\sg(a\dmin b)\lor\sg(b\dmin a))
\ \,\text{(de Morgan)}\\
&=\neg((a\dmin b)+(b\dmin a))\\
&=\neg(\dist(a,b))\\
&=\ttt{1}\dmin((a\dmin b)+(b\dmin a)):
\N\times\N\to\N\,.
\end{align*}

Substituting $\0$ for $b$ in the above deduction,
from the last line we get in particular
\begin{align*}
&(a\doteq\0)=\neg\,a:\N\to\N.
\end{align*}
\end{itemize}

\section{Order and equality}
\subsubsection*{Order and equality theorem}

\begin{enumerate}[(i)]
\item
Weak order $\leq:\N\times\N\to\two$ is
\begin{enumerate}[1.]
\item
\emph{reflexive}: $a\leq a$
\item
\emph{antisymmetric}: $(a\leq b)\land(b\leq a)
\impl(a\doteq b)$

\item
\emph{transitive}: $(a\leq b)\land(b\leq c)
\impl(a\leq c)$

and
\item
\emph{linear}: $(a\leq b)\lor(b\leq a)$
\end{enumerate}

\item
Strict order $<\,:\N\times\N\to\two$ is
\begin{enumerate}[1.]
\item
\emph{irreflexive}: $\neg\,(a<a)$\ \,and
\item
\emph{transitive}: $(a<b)\land(b<c)\impl(a<c)$
\end{enumerate}

\item
Weak order and strict order are connected by

$(a\leq b)=(a<b)\lor(a\doteq b)$

\item
The (binary) maximum is defined by
\begin{align*}
&\max(a,b)=a+(b\dmin a)=b+(a\dmin b)=\max(b,a)
\end{align*}

and is the (binary) supremum with respect to the order:
\begin{align*}
&(a\leq\max(a,b))\land(b\leq\max(a,b))\\
&\land
((a\leq c)\land(b\leq c)\impl(\max(a,b)\leq c))
\end{align*}

\item
The equality predicate `$\doteq$' is defined by 

$(a\doteq b)=(a\leq b)\land(b\leq a)$

\item
The equality predicate is an \emph{equivalence} predicate

i.e.\ it is 
\begin{enumerate}[1.]
\item
\emph{reflexive}:
$(a\doteq a)=(a\leq a)\land(a\leq a)=\true_\N$
\item
\emph{symmetric}: 
$(a\doteq b)\biimpl(b\doteq a)$\,,\ \,and
\item
\emph{transitive}: $(a\doteq b)\land(b\doteq c)
\impl(a\doteq c)$
\item
the triple 
$\set{<,\ \doteq,\ >}$ fulfils the 
\emph{law of trichotomy} 

$((a<b)\,\dot{\lor}\,(a\doteq b))\,\dot{\lor}\,(b<a)$

\emph{aut}(\emph{aut} $a<b$ \emph{aut} $a\doteq b$) 
\emph{aut} $b<a$
\end{enumerate}
\end{enumerate}

The proof is straight forward from the above,
in particular by use of truth tables, see section 4.4 above on Junctors and truth tables.

We now have at our disposition all ingredients for 
proving the

\section{Equality convertibility}
\subsubsection*{Equality definability theorem for $\PR$}

\inference{(\mr{EqDef})}
{$f=f(a):A\to B,\ g =g(a):A\to B$ in $\PR$\\
&$\PR\derives (f(a)\doteq_B g(a))$\,, $a$ free in $A$
}
{$\PR\derives\,f=g\ 
\ \text{\ie}\ \ f=^{\PR}g$}

In other words: A map equation which holds true predicatively 
for each argument
individually gives rise to an argument-free  
categorical equation
between the maps concerned.

\paragraph*{Proof.}
We begin with the special case 
$B = \N,\ \doteq_\N\ =\ (\,\doteq\,:\N\times\N\to\N)$\,. 

Let be $f,\ g: A \to \N$
$\PR$-maps satisfying the \emph{antecedent}
\begin{align*}
&\derives (f(a)\doteq g(a))
:A\to\N\times\N\xto{\doteq}\N\,, \,\text{hence}:\\
&\derives (g(a)\leq f(a))\land(f(a)\leq g(a))\\
&\derives \neg(g(a)\dmin f(a))
\ \text{and}\ \derives \neg(f(a)\dmin g(a))\\
&\derives(g(a)\dmin f(a))=\false_A(a)=\0_A(a)\ \text{and}\\ 
&\derives(f(a)\dmin g(a))
=\false_A(a)=\0_A(a)\\
& \text{therefore}\\
&\derives(g(a)\dmin f(a))=\0_A(a)\ \text{and}\ 
\derives(f(a)\dmin g(a))=\0_A(a)
\end{align*}

Replacing 
$(g(a)\dmin f(a))$ by $\0_A(a)$
resp.\ $\0_A(a)$ by $(f(a)\dmin g(a))$
we get
\begin{align*} 
\PR \derives
f&=f(a)+\0_A(a)\\
&=f(a)+(g(a)\dmin f(a))\\
&=\max(f(a),g(a))
\ \text{by {definition} of}\ \max(a,b) \\ 
&=\max(g(a),f(a))\ \text{by}\ \max
\ \text{symmetry} \\
&=g+(f(a)\dmin g(a))\\
&=g(a)+\0_A(a)=g(a)\\
&=g:A\to\N
\end{align*}

The general case where $B$ is a cartesian product of 
$\N$'s follows since individual equality on cartesian products is defined componentwise
\ \,\textbf{\qed}

By reflexivity of `$\doteq$' the converse inference to
$(\mr{EqDef})$ also 
to be true is shown in
\NAME{Pfender, Kröplin, Pape} 1994.

The conjunction of the inference $(\mr{EqConv})$ and its
converse inference gives the

\subsubsection*{Equality convertibility theorem}
\inference{(\mr{EqConv})}
{$f:A\to B$\,, $g:A\to B$ maps}
{$\PR\derives f=g$  
\textbf{iff}  $\PR\derives f(a)\doteq_B g(a)$\,, $a\in A$ free}

\subsubsection*{Proposition}
(Commutative) addition is strictly 
monotonous in both arguments. (Commutative) 
multiplication is (weakly monotonuos in both arguments
and) strictly monotonous in one argument provided the
other argument is greater than zero.
   
\paragraph*{Proof}
\begin{align*}
\begin{split}
(m<n)&\impl(\s m\leq n)\\
&\impl(\s m\dmin n\doteq\0)\\ 
&\impl((a+\s m)\dmin(a+n)\doteq\0)\\ 
&\impl(a+m<a+n)\ \text{as well as}
\end{split}\\
\begin{split}
(m<n)&\impl(\s m\leq n)\\
&\impl(\s m\dmin n\doteq\0)\\ 
&\impl((a\mul\s m)\dmin(a\mul n)\doteq 
a\mul(\s m\dmin n)\doteq\0)\\ 
&\impl(a\mul m<a\mul n)
\end{split}
\end{align*}

$\textbf{\qed}$

\section{Synopsis. Order from Algebra}

Alltogether we have 
\begin{enumerate}[(i)]
\item $\N$ admits the structure
\begin{align*}
\xymatrix @+1em{
\one
\ar @<+0.5ex> [r]^\0
\ar @<-0.5ex> [r]_{\s\0} 
& \N
  \ar[d]^{\sg}_{\,> \0}
&\N\times\N
\ar@<-0.6ex>@/_1pc/[l]_{\dmin}
\ar@<-0.5ex>[l]_{+}
\ar@<+0.5ex>[l]^{\mul}
\ar@<-0.5ex>[ld]_{\leq}
\ar@<+0.5ex>[ld]^{\doteq}
\\
&\two
}
\end{align*}

of a \emph{unitary commutative semiring with zero,} 
combined with

\item an additional algebraic
operation \emph{truncated subtraction} 
$m\dmin n:\N\times\N\to\N$ with its 
\emph{simplification property} $(a+n)\dmin(b+n)=a\dmin b$\,, 
and such that the multiplication \emph{distributes} over 
this kind of subtraction.

\item order
$(m\leq n):\N\times\N\to\two$ is reflexive, antisymmetric, 
transitive, and linear.

\item
$\max(a,b)=a+(b\dmin a)=b+(a\dmin b)=\max(b,a):
\N\times\N\to\N$

is the maximum with respect to the order.

\item
$\min(a,b)=\max(a,b)\dmin\abs{a-b}
=\max(b,a)\dmin\abs{b-a}$

$=\min(b,a):\N\times\N\to\N$ 

is the minimum with respect to the order.

Furthermore we have
\item
the fundamental equality \emph{predicate} 
\begin{align*}
&(m\doteq n)=\neg\,|m-n|:\N\times\N\to\two
\end{align*}

is an \emph{equivalence predicate} 
and makes up the
\emph{trichotomy}
 
$((a<b)\,\dot{\lor}\,(a\doteq b))\,\dot{\lor}\,(b<a)$\,, 
with strict order $(m<n)=(\s m\leq n)$

\item
Equality $f=g:A\to B$ \emph{between} maps is 
convertible (in $\PR$) with the equality \emph{predicate}
$(f(a)\doteq g(a)):A\times A\to\two$

\item Algebra Combined with Order:
Addition is strongly monotonous in both arguments, 
multiplication is strongly monotonous in one argument 
provided the other argument is greater than zero, 
and truncated subtraction is 
weakly monotonous in its first and weakly antitonous in 
its second argument.
\end{enumerate}

\section{Peano in $\bfS$}

We list the Peano axiom\index{axiom}s P1\,-\,P5 together with their 
\text{proofs} in $\bfS$\,.

\begin{itemize}
\item
P1: \emph{zero is a natural number}, $\0\in\N$: 

The map constant $\0:\one\to\N$ is a natural number.

\item
P2: \emph{To any natural number} $n$ 
\emph{is assigned a unique successor} $\s(n)$\,, 
$(\forall n\in\N) (\exists!\,n'\in\N)[n'=\s(n)]$\,:

As a p.r.\ map, $\s=\s(n):\N\to\N$ uniquely defines the 
successor of $n$\,.
In particular it defines a unique successor of any number 
constant.

\item
P3: $\0$ \emph{is not a successor}, 
$(\forall n\in\N)[\0\neq\s(n)]$:

We have $\s(n)>0$ and therefore $(\s(n)\neq\0):\N\to\N$\,,
\footnote{this is because the trichotomy (($(a<b)$
or else $(a\doteq b)$) or else $(a>b)$) holds in $\PR$
}

where $(a\neq b)=\neg(a\doteq b):\N\times\N\to\two$\,.
  
\item
P4: \emph{$m\neq n$ implies $\s(m)\neq\s(n)$}\,, or
by \emph{contraposition}, 

$(\forall m\in\N)
(\forall n\in\N)[\s(m)=\s(n)\impl m=n]$\,:

The successor map $\s:\N\to\N$ admits 
the predecessor map $\p:\N\to\N$ as a retraction:
$\p\circ\s=\id_\N$\,. Now since coretractions are 
injective, we have
$\bfS\derives(\s(m)\doteq\s(n))\impl(m\doteq n)$

\item
P5: Peano induction\index{induction} would require quantification over 
properties. 
In first-order logic this is circumvented by use
of an induction\index{induction} schema which corresponds to 
infinitely 
many \text{axiom\index{axiom}s}: 
\begin{align*}
&\ph(\0)\land(\forall n\in\N)[\ph(n)\impl\ph(n')]
\impl(\forall n\in\N)\ph(n)\,.
\end{align*}

This may also be written in the following form:

\inference{(\mr{P5})}
{$\ph=\ph(n):\N\to\two$\ predicate\\
&$\ph(\0)$\ (\emph{anchor})\\
&$(\forall n\in\N)[\ph(n)\impl\ph(\s(n))]$\ 
(\emph{induction\index{induction} step})
}
{$(\forall n\in\N)\ph(n)$ (\emph{conclusio})}

To this corresponds the Free-Variable $\PR$ schema

\inference{(\mr{P5_{FV\N}})}
{$\ph=\ph(n):\N\to\N$ predicate\\
&$\ph(\0)=\true:\one\to\N$\ (\emph{anchor})\\
&$(\ph(n)\impl\ph(\s(n)))=\true_\N(n):\N\to\N$\\
& (\emph{induction\index{induction} step})
}
{$\ph(n)=\true_\N(n):\N\to\N$ (\emph{conclusio})}

\paragraph*{Proof} of $(\mr{P5_{FV\N}})$\footnote{
  cf.\ \NAME{Pfender, Kröplin, Pape} 1994
  }:
 
Within schema $(\mr{pr}!)$ specify
\begin{align*} 
& g = g(a) = \ph(a,0) = \true_A(a):A\to\two\ (\emph{anchor})\\ 
& \text{and}\\
& h = h((a,n),b) = b\lor\ph(a,\s\,n)
:(A\times\N)\times\N\to\two\\
&  \ (\emph{induction\index{induction} step})
\end{align*}

By $(\mr{pr})$ together with (pr!) we obtain a unique 

  $f = f(a,n):A\times\N\to\two$ which satisfies
\begin{align*} 
& f(a,0)= \ph(a,0) = \true_A(a) \ \text{and} \\
& f(a,\s\,n) = h((a,n),f(a,n)) = f(a,n)\lor\ph(a,\s\,n) 
\end{align*}

This works for 
\begin{align*}
&f=\true_{A\times\N}:A\times\N\to\two\ \text{as well as for}\\
&f = \ph:A\times\N\to\two
\end{align*} 

the latter since
\begin{align*}
&\ph(a,n)\lor\ph(a,\s n)\\
&=(\ph(a,n)\lor\ph(a,\s\,n))\land(\ph(a,n)\impl\ph(a,\s n))\\
& \text{by antecedence (induction\index{induction} hypothesis)}\\
&=\ph(a,\s n)
\end{align*}

the last line by boolean tautology
\begin{align*}
&(\alpha\lor\beta)\land(\alpha\impl\beta)=\beta\ \text{(test with 
$\beta = \false$ and  $\beta=\true$)}
\end{align*} 

\textbf{\qed}

Making use of introducing the object $\two$ in the
$\PR$-extension $\bfS$ we may write the schema in the
more familiar boolean form
 
\inference{(\mr{P5_{FV}})}
{$\ph=\ph(n):\N\to\two$ predicate\\
&$\ph(\0)=\true$\ (\emph{anchor})\\
&$(\ph(n)\impl\ph(\s(n)))=\true$
\ (\emph{induction\index{induction} step})
}
{$\ph(n)=\true$ (\emph{conclusio})}

Furthermore, by replacing here the predicate $\ph$ by the 
conjunction
\begin{align*} 
&\underset{i \leq n}{\land}\ph(i):\N\to\two
\end{align*}

we get the free-variable

Course of values\index{course of values} induction\index{induction} principle\index{principle}

\inference{(\mr{P5^*_{FV}})}
{$\ph=\ph(n):\N\to\two$ predicate\\
&$\ph(\0)=\true$\ \,(\emph{anchor}) \\
&$(\underset{i\leq n}{\land}\ph(i) 
\impl\ph(\s(n)))=\true$\ \,(\emph{induction\index{induction} step})  
}
{$\ph(n)=\true$\ \,(\emph{conclusio})}
\end{itemize}

\chapter{Partial maps: from $\bfS$ to $\hatS$}
Turning to a categorical theory of partial maps we quote 
\NAME{Robinson, Rosolini} 1988 on partial functions:
``If one looks in a standard work on category theory, such as 
the classic \emph{Categories for the Working Mathematician}
(\NAME{MacLane} 1972),
one is hard put to it to find
examples of categories where the maps represent partial 
functions.'' 

\section{Recursive relations}
A \emph{(recursive) relation} is defined here as a
\emph{correspondence} (in $\bfS$) in the sense of
\NAME{Brinkmann, Puppe} 1969 i.e.\ as a \emph{span}\footnote{cf.\ \NAME{Robinson, Rosolini} 1988 
  and \NAME{Johnson, Rosebrugh} 2015 
  }, 
of form
\begin{align*}
&\xymatrix{
&R
\ar[dl]_l
\ar[dr]^r
\\
A
& & B
}
\end{align*}

This is formally different from the characterisation
of a \emph{non-empty relation} $f:W^r\to W^s$ in 
\NAME{Eilenberg, Elgot} 1970, \textsc{Corollary} 3.3 to their
\textbf{Main Theorem}: They use \emph{power sets} 
$\widehat{X}, \widehat{Y}$ not available in our 
\emph{equational} approach as well as ``arbitrary" unions
of sets.

In order to define a partial recursive \emph{function} 
we make use
of relations and of the definition of $\PR$\,, by closing 
under a suitable $\mu$-operation 
which is applied to a function rather than to a relation:

Suppose we are given a partial function $f(a,y)$\,. We 
introduce terms of form
\begin{align*}
&\mu(f(a,y))=
\begin{cases}
z\ \text{if $z$ in $\N$ is such that}\\
   f(a,z)\doteq\0\ \text{and for}\\
   w<z\ f(a,w)\neq0,\ w\ \text{free in $\N$}\,,\\
\text{undefined otherwise}
\end{cases}
\end{align*}

\section{Partial recursive maps}
A \emph{partial $\bfS$-map} or \emph{partial}
\emph{$\bfS$-function} $f:A\bs{\parto} B$ from \emph{domain}
$A$ to \emph{codomain} $B$ is to come here with an 
$\bfS$-object $D_f=\set{D:\delta_f}$\,, $D$ a 
$\PR$-object, i.e.\ a cartesian product of $\N$'s and
$\one$'s.

$D_f$ is mapped into domain $A$ by an $\bfS$-map
$d_f:D_f\to A$ which \emph{enumerates}, in general not injectively, the arguments for which $f$ is defined.

The \emph{values} of $f$ are specified
by an $\bfS$-map $\widehat{f}:D_f\to B$\,.

Thus we have the \textsc{diagram}
\begin{align*}
&\xymatrix @+0em{
D_f
\ar[d]^{d_f}
\ar[drr]^{\widehat{f}}
\\
A
&& B
} 
\end{align*}
which we call the
\emph{span} 
of the partial function $f:A\parto B$\,,
also written as the pair
\begin{align*}
&f=\an{d_f: D_f \to A,\ \widehat{f}: D_f \to B}
\end{align*}
or as the induced map\index{induced map}
\begin{align*}
&f=\an{(d_f,\widehat{f}):D_f\to A\times B}
\end{align*}

into the product $A\times B$\,.
  
The value $b=f(a)$ in $B$ of a defined argument $a$ in $A$
is to be $f(a)=\widehat{f}(\hat{a})$ where $a=d_f(\hat{a})$\,. 
 
To $a$ in $A$ correspond in general several $\hat{a}$ in $D_f$
with $d_f(\hat{a})=a$\,.
For \emph{well-definedness as a map}, 
i.e.\ \emph{functionality},

$f = \an{(d_f,\,\widehat{f})}$ is to fulfill, as part of 
its definition, the 

\emph{right-uniqueness} condition
\begin{align*}
&[d_f(\hat{a})\doteq d_f(\hat{a}')]
\impl[\widehat{f}(\hat{a})\doteq\widehat{f}(\hat{a}')]
\end{align*}

That is:
If an argument $a$ in $A$ has several inverse images $\hat{a}$ in 
$D_f$\,,\\
then $\widehat{f}$ will throw all of these on the same 
value\index{value} $b$ in $B$\,.

We define $\hatS$ to be the category of $\bfS$-objects 
and $\bfS$-partial functions.

\section{Inclusion and equality}
\emph{Inclusion}\index{inclusion} $f\,\parinc\,g$
of \emph{partial maps} is shown in the 
\textsc{diagram}

\begin{align*}
&\xymatrix@+0em{
&D_{f} 
\ar@/_1.0pc/[ldd]_{d_{f}} 
\ar[d]^i 
\ar@/^1.0pc/[rdd]^{\widehat{f}} 
\ar@{}[ldd]|{=}  
\ar@{}[rdd]|{=} 
& 
\\ 
&D_g
\ar[ld]_{d_g} 
\ar[rd]^{\widehat{g}} 
&&(\,\parinc\,)
\\
A 
\ar@{..^{>}}[rr]^{f}_g 
&&B
}
\end{align*}

where $i$ may be not injective. Map $i:D_f\to D_g$ 
\emph{defines} the inclusion.
   
The maps $f$ and $g$ in the above diagram are \emph{equal}, if
in the $\bfS$-\textsc{diagram} below
the equality sign `$=$' in the left part
also holds for a suitable map $j:D_g\to D_f$

\textsc{diagram}

\begin{align*}
&\xymatrix@+0em{
D_f
\ar[d]_{d_{f}}
\ar[rrrd]^{\widehat{f}}
\ar@<-.5ex>@/^2pc/[dd]_(0.6)i
\\
A
\ar@{}[rr]|(.15)=
&
\ar@{}[rr]|(.25)=
&&B
\\
D_g
\ar[u]^{d_g}
\ar@<-.5ex>@/_2pc/[uu]_(0.6)j
\ar[urrr]_{\widehat{g}}
}
\end{align*}

As implied by the diagram, the \emph{equality} $f\pareq g$ is 
\emph{reflexive}, \emph{symmetric}, and \emph{transitive}. 

With respect to this equality the \emph{inclusion} relation $f\,\parinc\,g$
is reflexive, transitive, and 
\emph{antisymmetric}.

The extended category $\hatS$ arises from
the following constructions.

\section{Composition of partial maps}
In order to define the \emph{composition} 
$g\parcirc f:A\parto B\parto C$ of two
partial maps $f:A\parto B$ and $g:B\parto C$ in $\hatS$ 
consider first the \textsc{diagram}
\begin{align*}
&\xymatrix@+1.0em{   
D_f 
\ar [d]^{d_f} 
\ar [rd]^{\widehat{f}}  
& D_g 
  \ar [d]^{d_g} 
  \ar [rd]^{\widehat{g}} 
  & 
\\
A 
\ar @{..^>} [r]^f 
& B 
  \ar @{..^>} [r]^g
  & C 
}
\end{align*}

which consists of two consecutive span\index{span}s.
  
In putting another span\index{span} on top of them we get \textsc{diagram} 

\begin{align*}
&\xymatrix@+0em{   
D_h 
\ar @/_2pc/ [dd]_{d_h}^{\ =} 
\ar [d]^{\pi_l} 
\ar [rd]^{\pi_r} 
\ar @/^2pc/ [rrdd]^{\widehat{h}}  
  & & 
\\
D_f 
\ar [d]^{d_f} 
\ar [rd]^{\widehat{f}} 
\ar @{} [r] | {\pb} 
& D_g 
  \ar [d]^{d_g} 
  \ar [rd]^{\widehat{g}} 
  \ar @{}[r] |(0.3) {=} 
  & 
\\
A 
\ar @{..^>} [r]^f 
\ar @/_1.5pc/ @{..^>} [rr]_{h\,=\,g\,\widehat{\circ} f}
& B 
  \ar @{..^>} [r]^g
  & C 
\\
&
}
\end{align*}

where the diamond $\mr{p.b.}$ is a 
\emph{pullback}.\footnote{
  \NAME{Robinson, Rosolini} 1988 likewise introduce partial
  maps as spans, and their composition via pullback. But
  both legs of their spans are to be monomorphisms. Our
  approach uses right-unique (``functional'') spans
  having arbitrary left legs. This is essential for 
  integrating $\mu$-recursive
  functions as partial functions
  }  

Composition $h=g\parcirc f:A \parto B\parto C$ 
is an $\hatS$-map.
\paragraph*{Proof:}
$h$ is right unique, since
\begin{align*}
&d_h(\hat{a})\doteq d_h(\hat{a}')\\
&\biimpl d_f(\pi_l(\hat{a}))
\doteq d_f(\pi_l(\hat{a}'))\\
&\impl\widehat{f}(\pi_l(\hat{a}))
\doteq\widehat{f}(\pi_l(\hat{a}'))
\quad(f\ \text{right unique})\\ 
&\biimpl d_g(\pi_r(\hat{a}))
\doteq d_g(\pi_r(\hat{a}'))\\
&\impl\widehat{g}(\pi_r(\hat{a}))
\doteq\widehat{g}(\pi_r(\hat{a}'))
\quad(g\ \text{right unique})\\
&\biimpl\widehat{h}(\hat{a})\doteq\widehat{h}(\hat{a}')
\ \textbf{\qed}
\end{align*}

\section{Monoidal product of partial maps}

\NAME{Budach, Hoehnke} 1975 
introduced the concept of \emph{half-terminal diagonal monoidal\footnote{
the main reference for \emph{monoidal products} is 
\NAME{Eilenberg, Kelly} 1966, where bicategories are used to generalise the concept of a
cartesian product to that of a \emph{tensor product}
`$\tensor$' which is \emph{left adjoint} to the 
\emph{internal \emph{hom} functor} and which in 
particular lacks in general the universal\index{universal} property of the 
cartesian product
} 
 category}
for the ``monoidal'' category of sets with partial 
maps, and denoted the product `$\tensor$'. 

Here this product of the p.r.\ 
partial maps $f:A\parto B$ and $f':A'\parto B'$ 
in $\hatS$ is written
$f\times f':A\times A'\parto B\times B'$ and 
defined by the span
\begin{align*}
&\xymatrix@+1em{
D_f\times D_{f'}
\ar[d]_{d_f\times d_{f'}}
\ar[dr]^{\widehat{f}\times\widehat{f'}}
\\
A\times A'
& B\times B'
}
\end{align*}

where $f\times f':A\times A'\parto B\times B'$ is a partial 
function.

\section{{Partial} {partial} maps}
For reduction of \emph{partial} partial p.\,r.\ maps to (just) 
partial p.\,r.\ maps consider \textsc{diagram}

\begin{minipage} {\textwidth}
$$
\xymatrix@+0em{
& D_{h} 
  \ar @{=} [d]^{\mr{def}}
  \ar @/_2pc/ [ddd]_(0.2){h=\gamma{h}}
  &
\\
& D_{\gamma{f}} 
  \ar @/_2 pc/ [lddd]_{d_{h} :\,= \,\widehat{d_f}}
  \ar [d]^{d_{\gamma{f}}}
  \ar @/^2pc/ [dd]^{\widehat{\gamma{f}}}
  \ar @/^2pc/ [rddd]^{\widehat{\widehat{f}}\, =\,: \widehat{h}}
  &
\\
& D_f          
  \ar @{-^{>}} @/^/ [u]^{d^-_{\gamma{f}}}
  \ar @{-^{>}} @/_/ [ldd]_{d_f}
  \ar @{-^{>}} [d]^{\gamma{f}}
  \ar @{-^{>}} @/^/ [rdd]^{\widehat{f}} 
  &
\\
& A \times B 
  \ar [ld]_{\ell}
  \ar [rd]^{\err}
  &
\\
A 
\ar @<0.5ex> @{-_{>}} [rr]^{f}
\ar @<-0.5ex> @{..^{>}} @/_1pc/ [rr]^{\ppareq}_{h}
  & {}
    & B
}
$$
\begin{center} Closure \textsc{diagram} for extension by partial maps 
  \end{center} 
\end{minipage}

\medskip   
The diagram shows a \emph{partial} partial p.\,r.\ map 
\begin{align*}
& f=\an{\gamma f:D_f\parto A\times B}:A \pparto B\\
& \text{-- arrow `$\pparto$'}
\ \text{for \emph{partial} partial --} 
\end{align*}

defined by its (partial) graph $\gamma f:D_f\parto A\times B$ 
in turn defined as an $\bfS$-map
\begin{align*}
& \gamma f = (d_{\gamma f},\widehat{\gamma f}):D_{\gamma f}\to A\times B
\end{align*}

As p.\,r.\ partial representant of \emph{partial} p.\,r.\ partial map 
$f: A \pparto B$ take the $\hatS$-map $h: A \parto B$ given by the \emph{frame} in the \textsc{diagram} above: 
\begin{align*}
h & = \an{(d_h,\widehat{h}): D_h \to A \times B}: A \parto B \\
& \defeq \an{(\widehat{d_f},\widehat{\widehat{f}}): 
                  D_{\gamma{f}} \to A \times B}: A \parto B \\
& \bydefeq \an{\widehat{\gamma{f}}: 
                  D_{\gamma{f}} \to A \times B}: A \parto B
\end{align*}

This shows: \emph{Partial} partial p.\,r.\ maps are 
(represented by) partial p.\,r.\ maps.

\section{The monoidal structure of $\hatS$}

\begin{itemize}
\item[$\mul$]
$\hatS$ is endowed with the following 
\emph{natural transformations}: 
\begin{align*}
&\mr{ASS}=(\ell\circ\ell,(\err\circ \ell,\err)):
(A \times B)\times C\overset{\iso}
{\lto}A\times(B\times C)\\
&\quad\emph{association}\\
&\Theta=(\err,\ell): 
A\times B\overset{\iso}{\lto}B\times A\ \emph{commutation}\\
&\Delta=(\id,\id):A\to A\times A\ \emph{diagonal}
\end{align*}

With these natural transformations
$\hatS$ can be shown to become a half-terminal\footnote{
  \emph{half}-terminal since the terminal map from $A$
  to $\one$ is not  
  uniquely determined} 
diagonal symmetric monoidal category:

\item[$\mul$]
\subsubsection*{Equality in monoidal category $\hatS$}
Composition and monoidal product of partial maps
are compatible with \emph{equality} 
`$\pareq$' of $\hatS$-maps.
\paragraph*{Proof}
Consider the \textsc{diagram}\footnote{
  we owe this diagram to F.\ Herrmann}
\begin{align*}
&
\xymatrix{ 
& & & D_{h'}
      \ar@{-->} '[lll] '[llldddd]^k [dddd]
      \ar[dr]|{\err'}
      \ar[dl]|{\ell'}
      \ar@/_1.5pc/[ddll]|{d_{h'}}
      \ar@/_-1.5pc/[ddrr]|{\widehat{h'}}
      \ar@{}[dd]|{\pb} 
\\
& & D_{f'}
    \ar@{-->}[dd]^(.2) {i}
    \ar[dr]|{\widehat{f}'}
    \ar[dl]|{d_{f'}} &
    & D_{g'}
      \ar@{-->}[dd]^(.2)j
      \ar[dr]|{\widehat{g}'}
      \ar[dl]|{d_{g'}} 
\\
& A
  \ar@{-^{>}}[rr]^(.7){f}_(.7){f'}
  \ar@{-^{>}}@/^1.2pc/[rrrr]^h|!{[dr];[rr]}\hole
                                        |!{[drrr];[rr]}\hole   
  \ar@{-^{>}}[rr]^(.7){f}_(.7){f'}
  \ar@{-^{>}}@/_1.2pc/[rrrr]_{h'}|!{[dr];[rr]}\hole
                                        |!{[drrr];[rr]}\hole   
  & & B
      \ar@{-^{>}}[rr]^(.3)g_(.3){g'}      
      & & C 
\\
& & D_{f}
    \ar[lu]|{d_f}
    \ar[ru]|{\widehat{f}} 
    & & D_{g}
        \ar[lu]|{d_g}
        \ar[ru]|{\widehat{g}} 
        & 
\\
& & & D_{h}
      \ar[ur]|{\err}
      \ar[ul]|{\ell}
      \ar@/_-1.5pc/[uull]|{d_{h}}
      \ar@/_1.5pc/[uurr]|{\widehat{h}}
      \ar@{}[uu]|{\pb}
}\\
&\text{Compatibility diagram 
                 \ of $\parcirc$ with $\subseteq$} 
\end{align*}

The partial map\index{partial map} $h:A\parto C$ is the 
composition\index{composition}
$\xymatrix{A\ar@{-^>}[r]^f&B\ar@{-^>}[r]^g&C}$

and the partial map $h':A\parto C$ is the composition\index{composition}

$\xymatrix{A\ar@{-^>}[r]_{f'}&B\ar@{-^>}[r]_{g'}&C}$

\item[$\mul$]
\subsubsection*{Associativity of composition
`$\parcirc$':} 
Consider the Associativity \textsc{diagram} 
for $\parcirc$ -- via \emph{nested pullbacks}

\begin{minipage} {\textwidth}
\begin{align*}
&\xymatrix{
D_{(h\widehat{\circ}g)\widehat{\circ}f}
\ar[dr]^{\pi_r}
\ar[d]^\iso_{\pi_l}
\ar@{}[ddr]|{\pb}
&&& 
\\
D_{h\widehat{\circ}(g\widehat{\circ}f)} 
\ar[d]_{\pi_l}
\ar[rrdd]^(.4){\pi_\err}|!{[rd];[rddd]}\hole 
&D_{h\widehat{\circ}g} 
\ar[dd]^(.45){\pi_\ell}
\ar[rdd]^{\pi_r} 
&& 
\\
D_{g\widehat{\circ}f}\ar[d]_{\pi_l} 
\ar[rd]\ar@{}[r]|{\pb} 
&&& 
\\
D_f 
\ar[d]^{d_f} 
\ar[rd]^{\widehat{f}} 
\ar @{} [r]|{\pb} 
&D_g 
\ar[d]^{d_g} 
\ar[rd]^{\widehat{g}}
\ar@{}[r]|{\pb}
&D_h 
\ar[d]^{d_h} \ar[rd]^{\widehat{h}} 
& 
\\
A 
\ar @{..^{>}} [r]^f & B 
\ar @{..^{>}} [r]^g & C 
\ar @{..^{>}} [r]^h & D
}
\end{align*}

  Associativity diagram for the composition\index{composition}
\end{minipage}

Here
$D_{(h\parcirc g)\parcirc f}
\xto{\pi_l}D_{h \parcirc(g\parcirc f)}$
is the pullback of
\begin{align*}
&\pi_r:D_{h\parcirc g}\to D_h\ \text{along}\ 
\pi_r:D_{h\parcirc(g\parcirc f)}\to D_h\,.
\end{align*}

This extends the \emph{cartesian} product of $\bfS$ 
into a \emph{bifunctor} on $\hatS$\,.
\end{itemize}

\subsubsection*{Inclusion Theorem}
Half-terminal diagonal monoidal category $\hatS$ admits a
corresponding inclusion functor $\mbf{I}:\bfS\lto\hatS$

\noindent
$\bfI$ is defined by
\begin{align*}
&\bfI(A\xto{f}B)
 =\xymatrix{
A
\ar[dr]^f
\ar[d]_\id
\\
A
&B
}
\end{align*}
The \textbf{proof} is straght forward.

\chapter{General Recursion}

\section{$\mu$-recursive functions within $\hatS$}

Given a function $f:A\times\N\to\N$ in $\bfS$\,, 
a $\mu$-recursive partial function $\mu f:A\parto\N$ within 
$\hatS$ is defined  by the span
\begin{align*}
&\xymatrix@+1em{
D_{\mu f}
\ar[dr]^{\widehat{\mu f}}
\ar[d]_{d_{\mu f}}
\\
A
\ar@{..^>}[r]^{\mu f}
&\N
}
\end{align*}

with components
\begin{align*}
&D_{\mu f}
=\set{(a,n)\in A\times\N:f(a,n)\doteq\0}\subseteq
A\times\N\,,\\
&d_{\mu f}=d_{\mu f}(a,n)=a
=\ell\,\circ\subseteq\,:D_{\mu f}\xto{\subseteq}
A\times\N\xto{\ell}A\,,\ \text{and}\\
&\widehat{\mu f}=\widehat{\mu f}(a,n)
=\min\set{m\leq n:f(a,m)\doteq\0}:D_{\mu f}
\xto{\subseteq}A\times\N\xto{\err}\N\,.
\end{align*}

That is: If an $a$ in $A$ is known to be 
coming with an $n$ 
in $\N$ such that $f(a,n)\doteq\0$\,, 
then $\mu f(a,n)$ is 
\emph{defined}, and $\mu f(a)$ is the \emph{minimal} 
$m\leq n$ with $f´(a,m)\doteq\0$\,.

Although for a given $a$ an $n$ 
such that $f(a,n)\doteq\0$ may not exist, \emph{functionality} 
of the span holds,
since for any $n, n'$ in $\N$ we have
\begin{align*}
&d_{\mu f}(a,n)\doteq d_{\mu f}(a,n')\\
&\impl
\widehat{\mu f}(a,n)
\doteq\widehat{\mu f}(a,\min(n,n'))
\doteq\widehat{\mu f}(a,n')
\end{align*}

i.e.\ the span is \emph{right-unique}.
Therefore $\mu$-recursive maps are in $\hatS$\,.

In case of a p.r.\ predicate $\ph=\ph(a,n):A\times\N\to\two$ we define 
a $\mu$-recursive partial function $\mu\ph:A\parto\N$ within 
$\hatS$ to be the span
\begin{align*}
&\xymatrix@+1em{
D_{\mu\ph}
\ar[dr]^{\widehat{\mu\ph}}
\ar[d]_{d_{\mu\ph}}
\\
A
\ar@{..^>}[r]^{\mu\ph}
&\N
}
\end{align*}

with components
\begin{align*}
&D_{\mu\ph}
=\set{(a,n)\in A\times\N:\ph(a,n)\doteq\true}\subseteq
A\times\N,\\
&d_{\mu\ph}=d_{\mu\ph}(a,n)=a
=\ell\,\circ\subseteq\,:D_{\mu\ph}\xto{\subseteq}
A\times\N\xto{\ell}A,\ \text{and}\\
&\widehat{\mu\ph}=\widehat{\mu\ph}(a,n)
=\min\set{m\leq n:\ph(a,m)\doteq\true}:D_{\mu\ph}
\xto{\subseteq}A\times\N\xto{\err}\N.
\end{align*}

That is: If an $a$ in $A$ is known to be 
coming with an $n$ 
in $\N$ such that $\ph(a,n)\doteq\true$\,, 
then $\mu\ph(a)$ is 
\emph{defined}, and $\mu\ph(a)$ is the \emph{minimal} 
$m\leq n$ with $\ph(a,m)\doteq\true$\,. 

Although for a given $a$ an $n$ 
such that $\ph(a,n)\doteq\true$ may not exist, \emph{functionality} 
of the span holds,
since for any $n, n'$ in $\N$ we have
\begin{align*}
&d_{\mu\ph}(a,n)\doteq d_{\mu\ph}(a,n')\\
&\impl
\widehat{\mu\ph}(a,n)
\doteq\widehat{\mu\ph}(a,\min(n,n'))
\doteq\widehat{\mu\ph}(a,n')
\end{align*}

i.e.\ the span is \emph{right-unique}.

\section{Partial $\bfS$-maps as $\mu$-recursive maps, $\hatS\bs{\iso}\muS$} 

\subsubsection*{Theorem}
To any (pointed\footnote{
  coming with a point map $\hat{a}_0:\one\to D_f$}) 
  \emph{partial} $\bfS$-map
\begin{align*} 
&f=\an{(d_f,\widehat{f}):D_f\to A\times B}:A\parto B
\end{align*} 

corresponds -- within $\hatS$ -- an
equal ($\pareq$) $\mu$-recursive $\hatS$-map
\begin{align*}
&g=\an{(d_g,\widehat{g}):D_g\to A\times B}:A\parto B
\end{align*}

\paragraph*{Proof.}
In order to specify such a $g$ we introduce the following
auxiliary p.r.\ maps:
\begin{align*}
&\enum_{D_f}:\N\to D_f\ \text{an enumeration},\\
&\ph_f=\ph_f(a,n):=(a\doteq_A\,d_f\circ\enum_{D_f}(n)): 
A\times\N\to\two\,,\\
&\mu\ph_f:A\times\N\parto\N\
\text{given by}\\
&D_{\mu\ph_f}=\set{A\times\N:\ph_f}
=\set{(a,n):a\doteq_A\,d_f\circ\enum_{D_f}(n)}
\subseteq A\times\N\\
&d_{\mu\ph_f}(\hat{a})=d_{\mu\ph_f}(a,n)=a=\ell\,\circ\subseteq
:\set{A\times\N:\ph_f}\xto{\subseteq}A\times\N\xto{\ell}A\\ 
& \text{and}\\
&\widehat{\mu\ph_f}=\widehat{\mu\ph_f}(a,n)=
\min\set{m\leq n:a\doteq_A\,d_{\mu\ph_f}(a,m)}\\
&=\min\set{m\leq n:a\doteq_A\,d_f\circ\enum_{D_f}(m)}\\
&:\set{A\times\N:\ph_f}\xto{\subseteq}A\times\N\xto{\err}\N
\end{align*}

We specify 
\begin{align*}
&g=\an{(d_g,\widehat{g}):D_g\to A\times B}:A\parto B
\end{align*}

as follows:
\begin{align*}
&D_g=D_{\mu\ph_f}\\
&d_g=d_{\mu\ph_f}:D_{\mu\ph_f}\to A\\
&\widehat{g}=\widehat{f}\circ\enum_{D_f}\circ\widehat{\mu\ph_f}
:D_{\mu\ph_f}\to\N\to D_f\to B
\end{align*}

In the \textsc{diagram}
\begin{align*}
&\xymatrix @+.0em{
&D_f
\ar@/_1pc/@<0.5ex>[ldd]^{i}
\ar[d]^{d_f}
\ar@<.5ex>[rrd]^{\widehat{f}}
\ar @{} [dr] |\pareq
&&
\\
{}
\ar@{}[r]|(0.65)=
&A
  \ar @<0.5ex> @{..^{>}} [rr]^(0.6){f}
  \ar @<-0.5ex> @{..^{>}} [rr]_(0.6){g}
  \ar @{-^{>}} [rd]^{\mu\ph_f}
  \ar @{} [d] |\pareq
  & & B
\\
D_{\mu\ph_f}
\ar @/^1pc/ @<0.5 ex> [ruu]^{j}
\ar [ru]|{\,d_{\mu\ph_f}}
\ar [rr]^{\widehat{\mu\ph_f}}
\ar@/_3pc/[rrru]_(.5){\widehat{g}}  
& & \N
      \ar [r]_{\enum}
      \ar @{}[ru]|\pareq
      & D_f
        \ar [u]_{\widehat{f}}
}
\end{align*} 

we take $i$ and $j$ to be
\begin{align*}
&i=i(a,n):D_f\to D_{\mu\ph_f}\\
&=(a,\min\set{m\leq n:
(d_f(a,m)\doteq_A d_f\circ\enum_{D_f}(m))}):\\
&D_f\to D_{\mu\ph_f}
=\set{(a,n):a\doteq_A\,d_f\circ\enum_{D_f}(n)}\\
& \text{and}\\ 
&j=j(a,n):=
\enum_{D_f}(\min\set{m\leq n:
d_f(a,m)\doteq_A d_f\circ\enum_{D_f}(m)}):\\
&D_{\mu\ph_f}\to D_f
\end{align*}

with the (dotted) partials $f$ and $g$\,.

\paragraph*{Corollary.} Category $\muS$ of $\mu$-recursion 
is monoidal, and it
is isomorphic to $\hatS$\,.
Both have cartesian p.r.\ category $\bfS$ included as a
half-terminal diagonal monoidal subcategory, and the inclusion 
is compatible with the isomorphism $\muS\bs{\iso}\hatS$\,.

\section{while-loop programs in $\hatS$}

A \ttt{while}-loop 
in $\hatS$ is a partial 
endomorphism $\wh[\ph|f]:A\parto A$\,, where 
$\ph=\ph(a):A\to\two$ is a p.r.\ \emph{control} 
predicate, and \emph{step} function $f:A\to A$ is
a p.r.\ endomorphism. 

Depending on the value $\ph(a)$ 
the iteration of $f$ terminates or continues. 

Before translating this into the span formalism, we may 
describe it in a conventional way: 
\begin{align*}
&\while[\ph|f](a) = 
\begin{cases}
a\ \textbf{if}\ \neg\ph(a)\,,\\
\while[\ph|f](f(a))\ \textbf{if}\ \ph(a)
\land\exists n\,(\neg\ph f^{n+\ttt{1}}(a))\,,\\
\emph{undefined}\ 
\textbf{if}\ \forall n\,(\ph f^{n+\ttt{1}}(a))
\end{cases}
\end{align*}

We convert this into the corresponding $\bfS$-span
\begin{align*}
&\xymatrix@+0em{
D_\wh
\ar[ddrr]^{\widehat{\mr{wh}}}
\ar[d]^\subseteq
\ar@/_2.2pc/[dd]_{d_\wh}^=
\\
A\times\N
\ar[d]^\ell
\\
A
\ar@{..^>}[rr]^\wh
&&A
}
\end{align*}

In linear form:
\begin{align*}
& \wh[\ph|f]=\an{(d_{\wh},\widehat{\mr{wh}}): 
D_{\wh}\to A\times A}:A\parto A\,,
\end{align*}

with the components
\begin{align*}
&D_\wh=\set{(a,n)\in A\times\N:\neg\ph\,f^\S(a,n)}
\subseteq A\times\N\\
&\quad\emph{defined arguments enumeration}\\
&d_\wh=d_\wh(a,n)
=\ell\circ\subseteq(a,n)=a:D_\wh\to A\\
&\quad\text{and \emph{rule}}\\
&\widehat{\wh}=\widehat{\wh}(a,n)
=f^\S(a,\min\set{m\leq n:
\neg\ph\,f^\S(a,m)}):D_\wh\to A
\end{align*}

The span defines an $\bfS$-partial \emph{endomap}.

\paragraph*{Proof.} It is to be shown that the 
span has the functional property of right-uniqueness:
\begin{align*}
&(a,n), (a',n')\in D_\wh=\set{(a,n)\in A\times\N:
\neg\ph\,f^\S(a,n)}\\
&\land(a\doteq d_\wh(a,n)\doteq d_\wh(a',n')\doteq a'\\
&\impl\widehat{\wh}(a,n)
\doteq\widehat{\wh}(a,\min(n,n'))\\
& \doteq\widehat{\wh}(a',\min(n,n'))
\doteq\widehat{\wh}(a',n')) 
\end{align*}

\textbf{\qed}  

\subsubsection*{Lemma. $\mbf{wh}$ by the $\bs{\mu}$-operator} 
\begin{align*}
&\wh[\chi|f]\pareq f^\S\parcirc 
(\id_A,\mu\ph_{\,[\chi|f]}):
A\parto A\times\N\to A\,,\ \text{where}\\
&\ph=\ph_{[\chi|f]}(a,n):=\neg\chi\circ f^\S(a,n): 
A\times\N\to A\to\two\,,
\end{align*}

therefore
\begin{align*}
&\wh[\chi|f]=f^\S\parcirc(\id_A,\mu[\chi\circ f^\S]):
A\parto A\times\N\to A\\
& \text{and, applied to an argument}\ a:\\
&\wh[\chi|f](a)=f^\S
\parcirc(\id_A,\mu\set{n:\chi\circ f^\S(a,n)}):
A\parto A
\end{align*}

\textbf{\qed}

\subsubsection*{Lemma. The $\bs{\mu}$-operator by $\mbf{wh}$}

For an $\bfS$-predicate $\ph=\ph(a,n):A\times\N\to\two$ 
the $\hatS$-function
\begin{align*}
&\mu\ph=\mu\ph(a):A\parto\N
\end{align*}

can be defined via $\wh$ as follows:
\begin{align*}
&\mu\ph:=\err
\parcirc\wh[\ph|A\times\s]\parcirc(\id_A,\0)\\
&:A\to A\times\N\parto A\times\N\xto{\err}\N\,.
\end{align*}

[As a computer program this definition of $\mu\ph$ reads
\begin{align*}
&\mu\ph(a):=\0\,;\\
&\mu\ph(a+\ttt{1})
:=\while\ \ph(a)\ \ttt{do}\ [a:=a+\ttt{1}]\ ]
\end{align*}

\textbf{\qed}

\section{Division with remainder by a while-loop}

We define 
here \emph{non-negative-integer division}
as a \ttt{while}-loop program by
\begin{align*}
&a\div b= 
\begin{cases}
n:=0\,;\\
\while\ a\leq n\ \ttt{do}\\
  [\ a:=a\dmin b\,;\\
  n:=n+\ttt{1}\ ];\\
\ttt{result}:=n
\end{cases}:\N\times\N\parto\N
\end{align*}

Formally:
\begin{align*}
&a\div b=
\err\parcirc\wh[a\leq n|\,
(a\dmin b,b)\times\s)]\parcirc((\N\times\N),\0)\\
&:\N\times\N\to\two\times\N\to\N\,,\\
& \text{with remainder}\\
&a\,\rem\,b=a\dmin(a\div b)\,.\\
& \text{It holds}\\
&a\doteq a\div b+a\,\rem\,b\,.
\end{align*}

This is obvious for $b\leq a$ but not trivial for $b>a$\,.

\section{Iterated functions in $\bfS$ as while-loops in $\hatS$}
 
An iterated $\bfS$-map 
$f^\S(a,n):A\times\N\to A$
is embedded into $\hatS$ 
as a \while-loop program:
\begin{align*}
&f^\S(a,n)= 
\begin{cases}
m:=0;\\
\while\ m\leq n\ \ttt{do}\\
 [\ a:=f(a);\\
  m:=m+\ttt{1}\ ];\\
\ttt{result}:=a
\end{cases}:A\times\N\to A
\end{align*}

Formally:
\begin{align*}
&f^\S(a,n)\pareq
\ell\parcirc\wh[m\leq n\,|\,(f\times\s)]
\parcirc(A\times\0)\\
&:A\times\N\to A\times\N\parto A\times\N\xto{\ell}A\,.
\end{align*}

This shows that any $\bfS$-iterated endomap has 
in $\hatS$ an equal \while-loop program.

\subsubsection*{Closure Remark}
By the above \emph{partial} partial maps equal in
$\hatS$ partial $\bf$ maps. Therefore the closure concerns
also \emph{nested} $\mu$-recursive maps
and \emph{nested} \while-loop programs.

\chapter{Total functions in $\hatS$}
Since $\hatS$ is an extension of $\bfS$ there are maps in
$\hatS$ which are total\index{total} functions. Because of the connectedness
of ``genuinely'' partial maps and total\index{total} functions the question
arises how to separate them from each other within $\hatS$\,.
The following Totality Lemma is to clarify this distinction.

\section{Totality Lemma}
\begin{enumerate}[]
\item[(i)]
For a partial p.r.\ map 
\begin{align*}
& f=\an{(d_f\,,\widehat{f}):D_f\to A\times B}:A\parto B
\end{align*}

the following statements are equivalent:
\begin{enumerate}[(a)]
\item
$f:A\parto B$ is -- in $\hatS$ -- equal to a \emph{total\index{total} map}

$g=\an{(\id_A,g):A\to A\times B}:A\to B$\,,

where $g$ is a suitable $\bfS$-map.
\item
$d_f=d_f(\hat{a}):D_f\to A$ is a retraction

\item
$d_f:D_f\to A$ admits the minimised opposite map
\begin{align*}
& d^-_f=d^-_f(a)
=\mu\set{\hat{a}\in D_f:d_f(\hat{a})\doteq a}:
A\to D_f
\end{align*}

as a -- total\index{total} -- coretraction. 
\end{enumerate}
\end{enumerate}

\subsubsection*{Proof.}
\begin{enumerate}[(i)]
\item[]
(a)$\to$(b): 

$f:A\parto B$ to be equal
$g$ is to say
\begin{align*} 
&\an{(d_f,\widehat{f})}\pareq\an{(\id_A,g)}:A\to B\,,
\end{align*}
  
which gives the span-inclusion 
\begin{align*} 
&\an{(\id_A,g)}\,\parinc\,\an{(d_f,\widehat{f})}\,.
\end{align*}

given by a map $i$
in the double-span diagram 
\begin{align*}
&\xymatrix@+0em{
&D_f
\ar[d]^{d_f}
\ar[dr]^{\widehat{f}}
\\
\ar@{}[r]|(.75)=
&A
\ar@{}[r]|=
&B&(*)
\\
&A
\ar[u]_(0.6){\id}
\ar[ur]_{g}
\ar@/^1.5pc/[uu]_(0.3)i
}
\end{align*}

This implies $d_f$ to be a retraction: assertion (b). 

(b)$\to$(c): Then the total\index{total} map
\begin{align*}
& d_f^- = d_f^-(a)
\bydefeq\mu\set{\hat{a}\in D_f:d_f(\hat{a})\doteq a}\\
&\pareq\min\set{\hat{a}\in D_f:\hat{a}\leq i(a)\,\land\,d_f(\hat{a})\doteq a}:
A\to D_f
\end{align*}

is equally a coretraction to $d_f:D_f\to A$\,, a
\emph{minimised} coretraction, 
constructed out of the coretraction $i:A\to D_f$
to $d_f$\,.
Thus (c) follows.
 
For (c)$\to$(a) consider the \textsc{diagram} below, 
with the ``suitable'' total\index{total} map
\begin{align*}
& \an{(\id_A,g)}
:=\an{(\id_A,\widehat{f}\circ d^-_f)}:A\to B\\
&\xymatrix@+0em{
\ar@{}[dr]|(.7)=
&A 
\ar@/_1.0pc/[ldd]_\id 
\ar@<-.5ex>[d]_{d_f^-}
\ar@/^1.0pc/[rdd]^{g=\widehat{f}\circ d_f^-}  
&
\\ 
&D_f
\ar@<-.5ex>[u]_{d_f} 
\ar[ld]_{d_f} 
\ar[rd]^{\widehat{f}} 
&&(**)
\\
A 
&&B
}
\end{align*}

We have 
$\an{(\id,g)}\,\parinc\,\an{(d_f,\widehat{f})}$
by coretraction property of $d_f^-$\,.

Conversely we have 
$\an{(d_f,\widehat{f})}\,\parinc\,\an{(\id,g)}$ 
because of
\begin{align*}
& g\circ d_f(\hat{a})=\\
&(\widehat{f}\circ d_f^-)\circ d_f(\hat{a})
=\widehat{f}(\hat{a}):D_f\to B, & (\bullet)\\
& \hat{a}\ \text{free in}\ D_f
\end{align*}

Because of right uniqueness of $\an{(d_f,\widehat{f})}$
we namely infer

\inference{}
{$d_f(\hat{a})=\id\circ d_f(\hat{a})$\\
&$=(d_f\circ d_f^-)\circ d_f(\hat{a})$\\
&$=d_f((d_f^-\circ d_f)(\hat{a})):D_f\to A$
}
{$\widehat{f}(\hat{a})
=\widehat{f}((d_f^-\circ d_f)(\hat{a}))$\\
&$=(\widehat{f}\circ d_f^-)\circ d_f(\hat{a}):D_f\to B$
}

Hence $(\bullet)$ holds and thus (c)$\to$(a).
Thus the equivalence of (a), (b), (c) is shown
\qed

\item[(ii)]
If in a composition
\begin{align*}
&h=g\parcirc f:A\parto B\parto C
\end{align*} 

the map $h$ is a total\index{total} map, then $f$ is a total\index{total} map.

\paragraph*{Proof.}
Consider 
the following commutative \textsc{diagram}.   
\begin{align*}
&\xymatrix@+.0em{   
A=D_{g\parcirc f} 
\ar[d]^{\pi_l} 
\ar[dr]^{\pi_r} 
\ar@<-0ex>@/_2pc/[dd]_(.6)\id
\ar@/^2.5pc/[ddrr]^{\widehat{h}}  
&& 
\\
D_f 
\ar[d]^{d_f} 
\ar[rd]^{\widehat{f}} 
\ar@{}[r]|{\mr{p.b.}} 
&D_g 
\ar[d]^{d_g} 
\ar[rd]^(.35){\widehat{g}} 
\\
A 
\ar@<-.3ex>@{..^>}[r]^f  
\ar@<-.7ex>@/_1pc/[rr]_{h\pareq g\parcirc f}
&B 
\ar@<-.3ex>@{..^>}[r]^g 
&C
}
\end{align*}

(The dotted arrows are to ease understending).
 
This diagram instantiates the earlier 
composition diagram, here for the total\index{total} map
$h=\an{(\id_A,h)}\pareq g\parcirc f$\,.

The $\bfS$-map $d_f$ is a retraction (to coretraction 
$\pi_l$).
Because of (b)$\to$(c)$\to$(a) this proves (ii). 

So in particular any coretraction in $\hatS$ is a total\index{total} map. 
\end{enumerate}
\textbf{\qed}

\section{Concluding Remarks}

Motivated by Dedekind and rejecting the unrestricted use of
quantifiers, Skolem\index{Skolem} 1923 instantiated what he called 
``recurrent way of thinking''. Goodstein\index{Goodstein} later realised Skolem\index{Skolem}'s aim to 
eliminate logical symbols completely in his own combinatorial
Foundation of Arithmetic, as likewise did Curry 
simultaneously. 

Here we construct the finitary Arithmetic of
Skolem\index{Skolem} and Goodstein\index{Goodstein}
categorically
without making use of concepts
which are not purely defined within category theory, not
presupposing set theoretical notions and not using logical
notation.

Aiming at a more general theory of recursive definability, 
Hilbert and
Bernays evoked Gödel's 1931 breakthrough,
still before Ackermann identified a recursive function
of a more general type than what was called
``Primitive Recursion'' since the work of P\'eter.

The functionals of $\PR$ are those of composition\footnote{
  here we adopt the notation of \emph{recursive functions} 
  in the Stanford Encyclopedia of Philosophy
  }
$\mc{C}omp=\mc{C}omp^{A,B}_C$
and
$\mc{P}rim\mc{R}ec=\mc{P}rim\mc{R}ec(\0,\mc{C}omp^{A,\N}_B)$\,.   

Categorically we get the
partial functions\index{partial functions} by the minimalisation\index{minimalisation} 
$\mu=\mu\ph:A\parto\N$
(operator``$\mc{M}in$'' in the Stanford Encyclopedia of
Philosophy),
transcending the realm of Primitive Recursion.

Here we get general recursion\index{general recursion} categorically
via partial functions
by the minimalisation
operator\index{minimalisation operator} ($\mc{M}in$) 
\begin{align*}
&\mu=\mu f:=
\begin{cases}
z\ \text{if}\ z\ \text{is such that}\ f(a,z)\doteq\0,\\
\text{for}\ w<z\ f(a,w)\neq\0,\ w\ \text{free},\\
\text{\emph{undefined} otherwise}.
\end{cases}
\end{align*}

This operation in combination with 
$\mc{C}omp=\mc{C}omp^{A,B}_C$ 
and $\mc{P}rim\mc{R}ec
=\mc{P}rim\mc{R}ec(\0,\mc{C}omp^{A,\N}_B)$
transcends the realm of Primitive Recursion which does
not have available $\mc{M}in$\footnote{
  see recursive functions in Stanford Encyclopedia of
  Philosophy, definitions 2.3, 2.4, 2.5}.
  
Nevertheless $\mu\bfS\bs{\iso}\hatS$ is a \emph{definitional} extension, containing
$\bfS$ embedded, and hence it is in particular a \emph{conservative}\index{conservative} extension of $\bfS$

\section*{Appendix. Free variables as projections}
\addcontentsline{toc}{section}
{Appendix. Free variables as projections}

We often make use of projections in order 
to get a \emph{variable-free}
notation. To enhance transparency of this proceeding
we describe it here in full detail, first by an example.

In a semiring\index{semiring} $R$ as e.g.\ $\N$ the distributive law\index{distributive law}
\begin{align*}
&a\mul(b+c)=(a\mul b)+(a\mul c)
\end{align*} 

corresponds to the commutative \textsc{diagram}
\begin{align*}
&
\xymatrix@+1em{
&(R\times R)\times R 
\ar@/_1pc/[dl]_(.7){(a,(b,c))}^(.7)
{(\ell\ell,(\err\ell,\err))}
\ar@/^1pc/ [dr]^(.7){((a,b),(a,c))}_(.7)
{((\ell\ell,\err\ell),(\ell\ell,\err))}
\ar[ddl]^(.7){(a,b+c)}
\ar[ddr]_(.7){(a \mul b,a \mul c)}
\\
R\times(R \times R)
\ar[d]_{R\,\times\,+}
&&(R\times R)\times(R \times R)
    \ar[d]^{\mul\,\times\,\mul}
\\
R\times R
\ar[dr]^{\mul}
&&R\times R
\ar[dl]_+
\\
&R
}
\end{align*}

wherin the free variable\index{variable}s are replaced by (composed)
projection\index{projection}s.

More generally this procedure can be described by
the following steps.
\begin{enumerate}[1.]
\item 
For the equation\index{equation}, extract the common codomain\index{common codomain} 
(domain of value\index{value}s), say $B,$  
(this codomain\index{codomain} may be implicit). 
\item 
Insert additional bracket pairs\index{bracket pairs} according to priority\index{priority} of 
operators.

\item 
Replace infix\index{infix} by prefix\index{prefix} notation.

\item 
Choose a linear order\index{linear order} for the occuring variable\index{variable}s
$a_1, \ldots, a_n$ as e.g,\ a lexical order\index{lexical order}. 

\item
Determine the common domain\index{common domain} for both sides of the equation\index{equation} -- to become the source of the diagram -- and write it in 
the form
\begin{align*}
&A=(\ldots((A_1\times A_2)\times\ldots)\times A_{\bs{m}})
\end{align*}

which in the context of this investigation is always
possible.
 
\item
\emph{Replace} all occurrence\index{occurrence}s
of a variable\index{variable} by the corresponding projection\index{projection}s.
 
\item
Use the projection operator \index{projection operator} $\Pi$ to replace numeral maps
$\0, \ttt{1}, \dots$ by $\0\,\Pi_A, \ttt{1}\,\Pi_A, \dots$
with the $A$ as above.
 
\item Insert the composition operator \index{composition symbol} 
`$\circ$' between terms which are not
bound together by an \emph{induced map\index{induced map} operator}
as in $(f_1,f_2)$\,.
\end{enumerate}

Also recall that for given
\begin{align*}
&a:=\ell_{A,B}:(A \times B)\to A\,,\\ 
&b:=\err_{A,B}:(A\times B)\to B\,\\
&f:A\to A',\ g: B \to B'
\end{align*}

the map product\index{map product} $f\times g$ is to be understood as
\begin{align*}
&(f\times g)(a,b)
=(f\times g)\circ(\ell_{A,B},\err_{A,B})\\
&=(f\times g)\circ\id_{A\times B}=f\times g\\
&=(f\circ\ell_{A,B},g\circ\err_{A,B})\\
&=(f\circ a,g\circ b)=(f(a),g(b))
:A\times B\to A'\times B'\,.
\end{align*} 

\chapter*{Outlook to Parts II and III}
\markboth{\ \hrulefill\ Outlook to Parts II and III}{\ Outlook to Parts II and III \hrulefill\ }
\addcontentsline{toc}{chapter}{Outlook to Parts II and III}

Work in preparation is Part II\index{Part II} on Evaluation using the 
categorical theory\index{categorical theory} $\hatS$ of 
\emph{partial maps\index{partial map}} to introduce 
general recursion\index{general recursion} and $\mu$-recursion, and the  
Totality Lemma\index{Totality Lemma} to separate totally defined
partial maps from general partial maps and to define 
evaluation\index{evaluation} by terminating 
\emph{Complexity Controlled Iteration\index{Complexity Controlled Iteration}},
under the roof of (categorical) Zermelo-Fraenkel\index{Zermelo-Fraenkel} \textbf{set} theory $\ZF$\,, for the latter see \NAME{Osius} 1974.

\noindent
Part III\index{Part III} in preparation is to deal with
\emph{Consistency Provability\index{Consistency Provability}}, to
discuss Gödel's\index{Gödel} second incompleteness theorem.  

\chapter*{Acknowledgements}
\markboth{\ \hrulefill\ Acknowledgements}{\ Acknowledgements \hrulefill\ }
\addcontentsline{toc}{chapter}{Acknowledgements}

My work on the present book dates back to the Vienna 2006
Conference celebrating Kurt Gödel's 100th anniversary.

Wolfgang Wenning helped me in a decisive way to improve this Part I of
the book.

Anton Ernst has examined the sections 
up to Primitive Recursion in this Part I.

Robert Sauer-Ernst organised the review of successive drafts by the AI 
assistant Claude (Anthropic) for editorial and consistency checking: 
cross-references, terminology, punctuation, and inconsistencies between 
chapters. All findings were assessed and, where accepted, implemented by 
the author, who takes full responsibility for the contents of the book, 
including suggestions adopted from his collaborators.

\chapter*{Index of Notation}
\markboth{\ \hrulefill\ Index of Notation}{\ Index of 
Notation\hrulefill\ }
\addcontentsline{toc}{chapter}{Index of Notation}

\begin{table}[htbp]
\begin{tabular*}{\textwidth}{@{\extracolsep{\fill}}l@{\enspace}l@{}}
\toprule
\multicolumn{2}{c}{\bf chapter 1 Cartesian language} \\
$\one$                  & terminal object one\\
$\N$                    & Natural Numbers Object NNO\\
$\times$                & cartesian product of sets and of maps\\
$\0$                & zero constant $0:\one\to\N$\\
$\s$                & successor function $\s:\N\to\N$\\
$\id$                   & identity map\\
$\circ$                 & map composition\\
$\Pi$                   & terminal map\\
$\ell/\err$             & left/right projection\\
$\mbf{var}$             & (free) variable\\
\bottomrule
\end{tabular*}
\end{table}

\begin{table}[htbp]
\begin{tabular*}{\textwidth}{@{\extracolsep{\fill}}l@{\enspace}l@{}}
\toprule
\multicolumn{2}{c}{\bf \textbf{chapter 2 Primitive Recursion}} \\
$\PR$                   & category/theory of primitive recursion\\
$\bfS=\PR+(\compr)$      & p.\,r.\ category/theory with predicate comprehension\\ 
                         & into subsets\\
$f^\S$                   & iteration of endo map $f$\\
$(\FR!)$                 & Freyd's uniqueness of the initialised iterated\\
$(\mr{pr})$              & full schema of primitive recursion\\
$[h*g]$             & p.\,r.\ map defined out of anchoring $g$\\
                         & and step map $h$\\
$\iso$                   & isomorphy\\ 
                         & resp.\ natural equivalence of functors\\
\bottomrule
\end{tabular*}
\end{table}

\begin{table}
\begin{tabular*}{\textwidth}{@{\extracolsep{\fill}}l@{\enspace}l@{}}
\toprule
\multicolumn{2}{c}{\bf \textbf{chapter 3 Goodstein Arithmetic GA}} \\                  
$\p$                   & predecessor map\\
$a\dmin b$               & truncated subtraction\\
$\mr{U}_1$ to $\mr{U}_4$ & Goodstein uniqueness rules\\
$\sg$                   & (arithmetical) sign\\
$\mr{V}_4$               & derived Goodstein uniqueness rule\\
$(a\doteq b)$    & individual equality, equality \emph{predicate}\\
$\neg$                   & negation\\
$\Delta$                 & diagonal map\\
\bottomrule
\end{tabular*}
\end{table}

\begin{table}
\begin{tabular*}{\textwidth}{@{\extracolsep{\fill}}l@{\enspace}l@{}}
\toprule
\multicolumn{2}{c}{\bf \textbf{chapter 4 Subobjects by predicates. Theory S}} \\                  
$\derives$               & (a theory) derives\\
$\p$                   & predecessor map\\
$\set{A:\chi}$           & subset abstracted from predicate $\chi$\\
$\bfS=\PR+(\compr)$      & p.\,r.\ theory $\PR$ + comprehension  subsets\\
$\two=\set{0,1}\subset\N$ & 2-element subset of NNO $\N$ \\
$f=^\bfS g$              & equality \emph{between} $\bfS$ maps \\
$\one+\N$                & sum/coproduct of objects\\
$A/\rho$                 & quotient set by an equivalence predicate\\
\bottomrule
\end{tabular*}
\end{table}

\begin{table}
\begin{tabular*}{\textwidth}{@{\extracolsep{\fill}}l@{\enspace}l@{}}
\toprule
\multicolumn{2}{c}{\bf \textbf{chapter 5 Order and equality on the NNO}} \\                  
$\alpha\sminus\beta$    & relative complement: $\alpha$ but not $\beta$\\
$\two=\set{\false,\true}$ & boolean algebra (in logical terms)\\
$\pret:\two\to\N$       & interpretation of truth values as numbers\\
$\sg:\N\to\two$       & logical signum of natural numbers\\
$\bfS=\PR+(\compr)$  & theory $\PR$ enriched\\ 
                          & by predicate comprehension into subsets\\
$f=\bs{if}[\chi,(h|g)](a):A\to B$ & definition of a map by case distinction\\
$\mr{P1}$ to $\mr{P5}$  & free-variable Peano axioms, here theorems,\\ 
                        & in particular\\
$\mr{P5}$               & Peano induction\\
\bottomrule
\end{tabular*}
\end{table}

\begin{table}
\begin{tabular*}{\textwidth}{@{\extracolsep{\fill}}l@{\enspace}l@{}}
\toprule
\multicolumn{2}{c}{\bf \textbf{chapter 6 Partial functions:
from $\bfS$ to $\hatS$}} \\                  
$\hatS$ &category of $\bfS$-objects\\ 
        &and $\bfS$-partial functions\\
$f=\an{d_f: D_f\to A,\ \widehat{f}:D_f\to B}$\\
$:A\parto B$&partial function\\
$[d_f(\hat{a})\doteq d_f(\hat{a}')]
\impl[\widehat{f}(\hat{a})\doteq\widehat{f}(\hat{a}')]$ &right-uniqueness condition\\
$f\,\parinc\,g$ & inclusion of partial maps\\
$f\pareq g$ & equality\\
$g\parcirc f:A\parto B\parto C$ & composition\\
$f\times f':A\times A'\parto B\times B'$ & monoidal product\\

$\mr{ASS}:(A \times B)\times C\overset{\iso}{\lto}A\times(B\times C)$
& association natural transformation\\
$\Theta=(\err,\ell):A\times B\overset{\iso}{\lto}B\times A$
& commutation transformation\\
$\Delta=(\id,\id):A\to A\times A$ & diagonal transformation\\
$\mbf{I}:\bfS\lto\hatS$ & inclusion functor\\

\bottomrule
\end{tabular*}
\end{table}

\begin{table}
\begin{tabular*}{\textwidth}{@{\extracolsep{\fill}}l@{\enspace}l@{}}
\toprule
\multicolumn{2}{c}{\bf \textbf{chapter 7 General Recursion}}\\                   
$\mu f:A\parto\N$ & $\mu$-recursive partial function\\
$\mu\bfS\iso\hatS$ & $\mu$-recursion $\iso$ partiality\\
$\wh[\ph|f]=\an{(d_{\wh},\widehat{\mr{wh}})}$ & while loop\\
$a\div b$ & integer division\\
$f^\S(a,n):A\times\N\to A$ & iterated $\bfS$-map\\
$\mu=\mu\ph:A\parto\N$ & minimalisation\\

\bottomrule
\end{tabular*}
\end{table}

\begin{table}
\begin{tabular*}{\textwidth}{@{\extracolsep{\fill}}l@{\enspace}l@{}}
\toprule
\multicolumn{2}{c}{\bf \textbf{chapter 8 Total functions in $\hatS$}}\\                  
$g=\an{(\id_A,g)}$ & total map\\
$d^-_f=d^-_f(a)$ & minimised opposite map\\
\bottomrule
\end{tabular*}
\end{table}

\printindex

\end{document}